\documentclass[11pt]{article}

\usepackage[margin=0.95in]{geometry}
\usepackage{mathtools}
\usepackage{amssymb}
\usepackage{amsthm}
\usepackage{dsfont}
\usepackage{stmaryrd}
\usepackage[T1]{fontenc}
\usepackage{microtype}

\usepackage{bm}
\usepackage[normalem]{ulem}
\makeatletter
\font\uwavefont=lasyb10 scaled 800
\renewcommand{\uwave}{\bgroup
  \UL@setULdepth
  \advance\ULdepth 4\p@
  \markoverwith{\lower\ULdepth\hbox{\uwavefont\char58}}\ULon}
\makeatother

\usepackage{enumitem}
\usepackage[nonamebreak]{natbib}
\usepackage[
    colorlinks=true,
    linkcolor=blue,
    citecolor=blue,
    urlcolor=blue
]{hyperref}
\usepackage{setspace}
\AtBeginDocument{
  \setlength{\abovedisplayskip}{7.5pt}
  \setlength{\belowdisplayskip}{7.5pt}
  \setlength{\abovedisplayshortskip}{7.5pt}
  \setlength{\belowdisplayshortskip}{7.5pt}
}
\newcommand\yestag{\addtocounter{equation}{1}\tag{\theequation}}

\newtheorem{theorem}{Theorem}[section]
\newtheorem{proposition}[theorem]{Proposition}
\newtheorem{corollary}[theorem]{Corollary}

\theoremstyle{definition}
\newtheorem{definition}[theorem]{Definition}
\newtheorem{remark}[theorem]{Remark}

\newcommand{\Pp}{\mathbb{P}}
\newcommand{\Ee}{\mathbb{E}}
\newcommand{\R}{\mathbb{R}}
\newcommand{\ind}{\perp\!\!\!\perp}
\newcommand{\cC}{\mathcal{C}}
\newcommand{\cE}{\mathcal{E}}

\newcommand{\cQ}{\mathcal{Q}}
\newcommand{\cU}{\mathcal{U}}
\newcommand{\cV}{\mathcal{V}}
\newcommand{\one}{\mathds{1}}

\makeatletter
\newcommand{\medl}[1]{\mathopen{\bBigg@{0.75}#1}}
\newcommand{\medr}[1]{\mathclose{\bBigg@{0.75}#1}}
\makeatother

\title{A Note on the Finite-Exceedance Dependence Principle for\\Sum--Maximum Asymptotic
Independence}
\author{Hongjian Shi}
\date{}

\begin{document}
\maketitle
\begin{abstract}
Asymptotic independence between sum-type and maximum-type statistics arises repeatedly across many
areas of probability and statistics, particularly in high-dimensional testing
and \linebreak random matrix theory. Motivated by the shared structure underlying existing proofs,
which is obscured by their model-specific formulations, this note isolates a common intermediate
principle.  We introduce a \emph{finite-exceedance dependence principle}:
if the \emph{finite-exceedance dependence error} vanishes at every fixed order
and high-order exceedances are mildly controlled, then a finite-order Bonferroni argument yields
asymptotic independence between the sum and the maximum.
We further give an independent-proxy criterion for verifying this:
for each finite set of exceedances, it suffices to construct an exactly independent proxy for the
sum-type statistic whose approximation error is negligible in aggregate.
Local deletion under dependency graphs and Gaussian regression and residualization are two examples.
The framework gives modular derivations of existing \linebreak sum--max independence results.
\end{abstract}


\section{Introduction}

The joint asymptotic behavior of sums and maxima is a classical topic in statistics.
The sum is typically \linebreak driven by the aggregate contribution of many moderate observations,
whereas the maximum is determined by a small number of extreme ones. It has long been known that,
for independent observations, this separation of scales leads to asymptotic independence under
suitable tail assumptions; see, e.g.,\linebreak
\citet*{Oliveira1962,Rosengard1962,Rossberg1965,ChowTeugels1979}.
Related results have also been developed for strongly mixing and weakly dependent stationary
sequences and stationary Gaussian sequences (see, among others,
\citet*{AndersonTurkman1991,AndersonTurkman1993,AndersonTurkman1995,Hsing1995a,Hsing1995b,HoHsing1996}),
and, in related forms, for both discrete and continuous-time Gaussian processes
(e.g., \citet*{HoMcCormick1999,McCormickQi2000}).

A closely related phenomenon has become increasingly important in high-dimensional statistics.
Sum-type statistics are powerful against dense alternatives, while max-type statistics are powerful
against sparse alternatives. Their asymptotic independence permits adaptive combination procedures.
Examples include high-dimensional independence testing \citep*{LiXue2015},
adaptive $U$-statistic tests \citep*{HeEtAl2021}, one-sample mean problem
\citep*{ChenFengLiang2024}, finite-$L_q$ rank blocks combined with maxima \citep*{ZhaoWangFeng2026},
and dependent Gaussian sum--max tests \citep*{FengEtAl2024}.
A point-process version of the same sum--maximum separation also appears in a recent work on
random matrix theory \citep*{HeinyKleemann2023}.
\clearpage

Although the final conclusions look similar and the proofs are often presented in problem-specific language,
a recurring structure can be identified as follows.

\begin{enumerate}[label=(\roman*),itemsep=-.1em, topsep=.5em]
\item Write the max event as a union of rare exceedance events
    and apply Bonferroni inclusion--exclusion.
\item For a finite set of candidate extreme coordinates, construct a modified sum-type statistic
    that is \emph{exactly independent} of those exceedance events.
\item Show that the modification is negligible on the scale of the sum-type statistic,
    both uncondi\-tionally and after regression on and summing over the rare exceedance
    configurations.
\item Use anti-concentration of the sum-type statistic to transfer the exact independence of the
    proxy back to the original statistic.
\end{enumerate}

The purpose of this note is to isolate this argument at what we regard as a useful level of
abstraction. We would like to characterize when and why the asymptotic factorization of every finite
exceedance intersection occurs. The balanced formulation is therefore an \emph{independent-proxy
criterion} combined with a generic finite-exceedance Bonferroni-type lemma.

The core result has two layers. First, define the aggregated \emph{finite-exceedance dependence
error}~as
\[
D_{n,k}(y)
=
\sup_{x\in\R}
\sum_{|I|=k}
\left|
\Pp\{T_n\le x,\,\cE_{n,I}(y)\}
-
\Pp\{T_n\le x\}\Pp\{\cE_{n,I}(y)\}
\right|,
\]
where $T_n$ is a real-valued statistic of interest (in this paper, most often a sum-type statistic)
and $\cE_{n,I}(y)$ denotes the intersection of $k$ exceedance events $\left\{X_{n,i}>u_n(y)\right\}$ indexed by~$i\in I$.
Theorem~\ref{thm:bonferroni} then uses $D_{n,k}(y)$ to provide a bound for
\[
\sup_{x\in\R}
\left|
\Pp\{T_n\le x,\cU_n(y)\}
-\Pp\{T_n\le x\}\Pp\{\cU_n(y)\}
\right|
\]
where $\cU_n(y)=\left\{\max_{i}X_{n,i}>u_n(y)\right\}$ is the exceedance event for the corresponding
maximum. Thus, showing that $D_{n,k}(y)\to0$ for each fixed $k$,
together with a mild high-order binomial-moment condition, is sufficient for sum--max asymptotic
independence. Second, Theorem~\ref{thm:proxy} gives a reusable sufficient condition for
$D_{n,k}(y)\to0$ based on an exactly independent proxy. Therefore,
the problem is reduced to constructing an appropriate independent proxy
and verifying the conditions of Theorem~\ref{thm:proxy}.

When a sparse dependency graph can be formulated, the local-deletion argument in
Section~\ref{sec:deletion} is a natural way to construct an appropriate proxy.
This mechanism is used in all examples in Section~\ref{sec:applications}
except Section~\ref{subsec:dpdGau}, where Gaussian regression and residualization provide an alternative construction.

\enlargethispage*{.5\baselineskip}

\section{Finite exceedances and the Bonferroni reduction}

Let $d_n\to\infty$ as $n\to\infty$ and let \(\cE_{n,1}(y),\ldots,\cE_{n,d_n}(y)\) be
exceedance events indexed by a fixed threshold parameter $y\in\R$. In applications,
one typically has
\begin{equation}\label{eq:EE}
\cE_{n,i}(y)=\{X_{n,i}>u_n(y)\},
\end{equation}
where random variables $X_{n,1},\ldots,X_{n,d_n}$ are not necessarily i.i.d.
The union of all exceedance events
\[
\cU_n(y):=\bigcup_{i=1}^{d_n}\cE_{n,i}(y)
\]
is exactly the event that the maximum $\max_{1\le i\le d_n}X_{n,i}$ exceeds the level $u_n(y)$
when the $\cE_{n,i}(y)$ are exceedance events given in \eqref{eq:EE}.

Let $T_n$ be an arbitrary real-valued statistic and let
\[
F_n(x):=\Pp(T_n\le x).
\]
For $I\subseteq\llbracket d_n\rrbracket:=\{1,\ldots,d_n\}$, define
\[
\cE_{n,I}(y)=\bigcap_{i\in I}\cE_{n,i}(y).
\]
For $k\ge1$, define the $k$th exceedance probability sum
\begin{equation*}
H_{n,k}(y)
=
\sum_{|I|=k}
\Pp\{\cE_{n,I}(y)\}
\equiv
\sum_{\substack{I\subseteq\llbracket d_n\rrbracket,\,|I|=k}}
\Pp\{\cE_{n,I}(y)\}.
\end{equation*}
Let $N_n(y)=\sum_{i=1}^{d_n}\one\{\cE_{n,i}(y)\}$ denote the exceedance count. Then $H_{n,k}(y)$ is the $k$-th binomial moment of $N_n(y)$ \citep*{Sibuya1991}
\[
H_{n,k}(y)=\Ee \binom{N_n(y)}k=\frac{1}{k!}\Ee\bigl[(N_n(y))_k\bigr],
\]
where $(a)_k:=a(a-1)\cdots(a-k+1)$ is the falling factorial.

The central quantity in this note is the following.

\begin{definition}[Finite-exceedance dependence error]
For $k\ge1$, define
\begin{equation*}
D_{n,k}(y)
=
\sup_{x\in\R}
\sum_{|I|=k}
\left|
\Pp\{T_n\le x,\cE_{n,I}(y)\}
-
\Pp\{T_n\le x\}\Pp\{\cE_{n,I}(y)\}
\right|.
\end{equation*}
\end{definition}

The following theorem establishes asymptotic independence
under $D_{n,k}(y)\to0$ for each fixed $k$,
together with a mild condition on the higher-order binomial moment $H_{n,2r}(y)$.
We defer discussion for proving $D_{n,k}(y)$ small to Theorem~\ref{thm:proxy}.

\begin{theorem}[Finite-exceedance dependence principle]\label{thm:bonferroni}
For every integer $r\ge1$ and every $y\in\R$,
\begin{equation}\label{eq:finite-bound}
\sup_{x\in\R}
\left|
\Pp\{T_n\le x,\cU_n(y)\}
-\Pp\{T_n\le x\}\Pp\{\cU_n(y)\}
\right|
\le
\sum_{k=1}^{2r}D_{n,k}(y)+H_{n,2r}(y).
\end{equation}
Consequently, if for every fixed $k\ge1$,
\begin{equation}\label{eq:D-to-zero}
D_{n,k}(y)\longrightarrow0,
\end{equation}
and
\begin{equation}\label{eq:H-tail}
\lim_{r\to\infty}\limsup_{n\to\infty}H_{n,2r}(y)=0,
\end{equation}
then
\begin{equation}\label{eq:asymp-factor}
\sup_{x\in\R}
\left|
\Pp\{T_n\le x,\cU_n(y)\}
-\Pp\{T_n\le x\}\Pp\{\cU_n(y)\}
\right|
\longrightarrow0.
\end{equation}
The same conclusion holds with $\cU_n(y)$ replaced by its complement.
\end{theorem}
\clearpage

\begin{proof}
For $k\ge1$, define
\[
J_{n,k}(x,y)
=
\sum_{\substack{I\subseteq\llbracket d_n\rrbracket\\|I|=k}}
\Pp\{T_n\le x,\cE_{n,I}(y)\}.
\]
Bonferroni's inequalities imply
\begin{align}
\sum_{k=1}^{2r}(-1)^{k-1}J_{n,k}(x,y)
&\le
\Pp\{T_n\le x,\cU_n(y)\}
\le
\sum_{k=1}^{2r-1}(-1)^{k-1}J_{n,k}(x,y),
\label{eq:bonf-joint}
\end{align}
and
\begin{align}
\sum_{k=1}^{2r}(-1)^{k-1}H_{n,k}(y)
&\le
\Pp\{\cU_n(y)\}
\le
\sum_{k=1}^{2r-1}(-1)^{k-1}H_{n,k}(y).
\label{eq:bonf-max}
\end{align}
The difference between the odd and even Bonferroni truncations in \eqref{eq:bonf-max} is exactly
$H_{n,2r}(y)$.

By definition of $D_{n,k}(y)$,
\[
\left|J_{n,k}(x,y)-F_n(x)H_{n,k}(y)\right|
\le D_{n,k}(y)\qquad\text{uniformly in $x$.}
\]
Combining this inequality with the upper bound in \eqref{eq:bonf-joint}
and the lower bound in \eqref{eq:bonf-max} gives
\begin{equation}\label{eq:finite-bound1}
\Pp\{T_n\le x,\cU_n(y)\}-\Pp\{T_n\le x\}\Pp\{\cU_n(y)\}
\le
\sum_{k=1}^{2r-1}D_{n,k}(y)+H_{n,2r}(y).
\end{equation}
Similarly, combining the lower bound in \eqref{eq:bonf-joint} and the upper bound in
\eqref{eq:bonf-max} yields
\begin{equation}\label{eq:finite-bound2}
\Pp\{T_n\le x\}\Pp\{\cU_n(y)\}-\Pp\{T_n\le x,\cU_n(y)\}
\le
\sum_{k=1}^{2r}D_{n,k}(y)+H_{n,2r}(y).
\end{equation}
Taking the supremum on both \eqref{eq:finite-bound1} and \eqref{eq:finite-bound2} over $x$ proves
\eqref{eq:finite-bound}.

For fixed $r$, assumption \eqref{eq:D-to-zero} makes the finite sum on the right-hand side vanish as
$n\to\infty$. Taking $\limsup_n$ in \eqref{eq:finite-bound} and then sending $r\to\infty$ proves
\eqref{eq:asymp-factor}. The statement for $\cU_n(y)^c$ follows from
\[
\Pp\{T_n\le x,\cU_n(y)^c\}-\Pp\{T_n\le x\}\Pp\{\cU_n(y)^c\}
=-\Bigl[
\Pp\{T_n\le x,\cU_n(y)\}-\Pp\{T_n\le x\}\Pp\{\cU_n(y)\}
\Bigr].
\]
\end{proof}

\begin{remark}[When a fixed Bonferroni order is sufficient—and
when it is not]\label{rem:fixed-order}
When the goal is to establish asymptotic independence, Theorem~\ref{thm:bonferroni} permits the
order $r$ to be fixed while $n\to\infty$, followed only afterwards by $r\to\infty$.
This is often substantially simpler than choosing a truncation order that grows with dimension.
However, it is worth emphasizing that the Bonferroni bound given in \eqref{eq:finite-bound} is
non-asymptotic. Choosing a growing order \(r\) may still be important for obtaining an explicit
uniform convergence rate, but it is not necessary to the asymptotic independence itself.
\end{remark}

\bigskip

Theorem~\ref{thm:bonferroni} establishes asymptotic independence
but does not, by itself, guarantee the existence of a joint limiting distribution.
The following corollary fills this gap by providing a more self-contained approach
to establishing both the marginal limit of the maximum
and its asymptotic independence from the sum.

\newpage
\begin{corollary}[Poisson criterion]\label{cor:poisson}
Suppose that for every fixed $k\ge1$ and $y\in\R$,
\begin{equation}\label{eq:poisson-fm}
H_{n,k}(y)\longrightarrow \left.{\{\lambda(y)\}^k}\right/{k!}
\end{equation}
for some $0\le\lambda(y)<\infty$, and that \eqref{eq:D-to-zero} holds. Then
\begin{equation}\label{eq:marginal-poisson}
\Pp\{\cU_n(y)^c\}\longrightarrow e^{-\lambda(y)}
\end{equation}
and
\begin{equation}\label{eq:joint-poisson}
\sup_{x\in\R}
\left|
\Pp\{T_n\le x,\cU_n(y)^c\}
-F_n(x)e^{-\lambda(y)}
\right|
\longrightarrow0.
\end{equation}
If in addition $T_n\Rightarrow T$ with continuous distribution function $F$,
then for every $x,y\in\R$,
\begin{equation}\label{eq:joint-poisson-ae}
\Pp\{T_n\le x,\cU_n(y)^c\}\longrightarrow F(x)e^{-\lambda(y)}.
\end{equation}
\end{corollary}

\begin{proof}
Equation \eqref{eq:poisson-fm} implies \eqref{eq:H-tail}.
The Bonferroni argument applied only to $\cU_n(y)$ yields $\Pp\{\cU_n(y)^c\}\to e^{-\lambda(y)}$,
while Theorem~\ref{thm:bonferroni} gives the factorization~\eqref{eq:joint-poisson}.
If $F_n\to F$ weakly and $F$ is continuous, then the convergence of distribution functions is
everywhere. This fact together with \eqref{eq:joint-poisson} concludes \eqref{eq:joint-poisson-ae}.
\end{proof}

\section{An independent-proxy criterion}\label{sec:proxy}

Theorem~\ref{thm:bonferroni} isolates the exact quantity to be controlled.
We now give a general mechanism for proving this vanishing.
This together with Theorem~\ref{thm:bonferroni} yields the main reusable criterion of the note.

\begin{theorem}[Independent-proxy criterion]\label{thm:proxy}
Fix $k\ge1$ and $y\in\R$. Assume that, for every $I\subseteq\llbracket d_n\rrbracket$ with $|I|=k$,
a proxy statistic $T_{n,I}^{\circ}$, which may depend on $I$ and need not have the same distribution
for different~$I$, can be constructed such that
\begin{equation*}
T_{n,I}^{\circ}\ind \cE_{n,I}(y).
\end{equation*}
For $\eta>0$, define
\begin{align*}
A_{n,k}(\eta,y)
&=
\sum_{|I|=k}
\Pp\{|T_n-T_{n,I}^{\circ}|>\eta\}
\,\Pp\{\cE_{n,I}(y)\},
\\
B_{n,k}(\eta,y)
&=
\sum_{|I|=k}
\Pp\{|T_n-T_{n,I}^{\circ}|>\eta,\cE_{n,I}(y)\},
\\
\intertext{and the Lévy concentration function}
\omega_n(\eta)
&=
\sup_{x\in\R}\{F_n(x+\eta)-F_n(x-\eta)\}.
\end{align*}
Then, for every $\eta>0$,
\begin{equation}\label{eq:proxy-bound}
D_{n,k}(y)
\ \le\ 
\underbrace{A_{n,k}(\eta,y)}
_{\substack{\text{weighted unconditional}\\
            \text{proxy error}}}
+
\underbrace{B_{n,k}(\eta,y)}
_{\substack{\text{weighted conditional}\\
            \text{proxy error}}}
+\quad
\underbrace{\omega_n(2\eta)\times H_{n,k}(y)}
_{\text{anti-concentration}}.
\end{equation}
Consequently, if there exists $\eta_n\to0$ such that
\begin{equation}\label{eq:proxy-conditions}
A_{n,k}(\eta_n,y)\to0,\qquad
B_{n,k}(\eta_n,y)\to0,\qquad
\omega_n(2\eta_n)\to0,\qquad
H_{n,k}(y)=O(1),
\end{equation}
then $D_{n,k}(y)\to0$.
\end{theorem}

\begin{proof}
Fix $I\subseteq\llbracket d_n\rrbracket$ with $|I|=k$. Abbreviate $\cE_I=\cE_{n,I}(y)$ and $T_I^{\circ}=T_{n,I}^{\circ}$.
Let
\[
\alpha_I(\eta)=\Pp(|T_n-T_I^{\circ}|>\eta),
\qquad
\beta_I(\eta)=\Pp(|T_n-T_I^{\circ}|>\eta,\cE_I),
\qquad
p_I=\Pp(\cE_I).
\]
It holds that
\begin{equation*}
T_n\le x,
\quad |T_n-T_I^{\circ}|\le\eta
\quad\Longrightarrow\quad
T_I^{\circ}\le x+\eta.
\end{equation*}
Hence, using this triangle inequality,
\begin{align*}
\Pp(T_n\le x,\cE_I)
&\le
\Pp(T_I^{\circ}\le x+\eta,\cE_I)+\Pp(|T_n-T_I^{\circ}|>\eta,\cE_I)\\
&=
 \Pp(T_I^{\circ}\le x+\eta)\,p_I+\beta_I(\eta),\yestag\label{eq:firstApprox}
\end{align*}
where the last step is by the exact independence $T_{n,I}^{\circ}\ind \cE_{n,I}(y)$.
Using the triangle inequality once again gives
\begin{align*}
\Pp(T_I^{\circ}\le x+\eta)
&\le \Pp(T_n\le x+2\eta)+
\Pp(|T_n-T_I^{\circ}|>\eta)\\
&=
 F_n(x+2\eta)+\alpha_I(\eta).\yestag\label{eq:secondApprox}
\end{align*}
Combining \eqref{eq:firstApprox} and \eqref{eq:secondApprox} shows
\begin{equation}\label{eq:proxy-upper}
\Pp(T_n\le x,\cE_I)-p_I F_n(x)
\le
\omega_n(2\eta)\,p_I+\alpha_I(\eta)\,p_I+\beta_I(\eta).
\end{equation}

\enlargethispage*{.5\baselineskip}

Similarly, the implication
\[
T_I^{\circ}\le x-\eta,
\quad |T_n-T_I^{\circ}|\le\eta
\quad\Longrightarrow\quad
T_n\le x
\]
yields
\[
\Pp(T_n\le x,\cE_I)
\ge
\Pp(T_I^{\circ}\le x-\eta)\,p_I-\beta_I(\eta),
\]
and
\[
\Pp(T_I^{\circ}\le x-\eta)
\ge
F_n(x-2\eta)-\alpha_I(\eta).
\]
Therefore
\begin{equation}\label{eq:proxy-lower}
p_I F_n(x)-\Pp(T_n\le x,\cE_I)
\le
\omega_n(2\eta)\,p_I+\alpha_I(\eta)\,p_I+\beta_I(\eta).
\end{equation}
Combining \eqref{eq:proxy-upper} and \eqref{eq:proxy-lower}, summing over $|I|=k$,
and then taking the supremum over $x$ proves \eqref{eq:proxy-bound}.
\end{proof}

\begin{remark}[The two approximation errors are genuinely different]\label{rem:two-errors}
The term $A_{n,k}$ is the \emph{weighted unconditional proxy error}:
it weights the unconditional probability that the proxy approximation fails by the probability of
the corresponding rare configuration. The term $B_{n,k}$ is the \emph{weighted conditional proxy
error}:
\[
B_{n,k}(\eta,y)
=
\sum_{|I|=k}
\Pp\{|T_n-T_{n,I}^{\circ}|>\eta \mid \cE_{n,I}(y)\}
\,\Pp\{\cE_{n,I}(y)\},
\]
which weights the conditional proxy-failure probability given the rare configuration by the
probability of that configuration. The latter is the quantity that a purely unconditional
$o_{\Pp}(1)$ statement can miss. \linebreak In many applications both terms are controlled by an
extremely strong uniform tail estimate, but the distinction is useful in more delicate problems.
\end{remark}
\clearpage

\begin{corollary}[Uniform proxy criterion]\label{cor:uniformproxy}
Define the following uniform bound for $T_n-T_{n,I}^{\circ}$
\begin{equation}\label{eq:rhonk}
\rho_{n,k}(\eta)
:=\sup_{|I|=k}\Pp(|T_n-T_{n,I}^{\circ}|>\eta).
\end{equation}
Then, for every $\eta>0$,
\begin{equation}\label{eq:D-bound}
D_{n,k}(y)
\le
\rho_{n,k}(\eta)\,H_{n,k}(y)
+\rho_{n,k}(\eta)\,\binom {d_n}k
+\omega_n(2\eta)\times H_{n,k}(y)
\end{equation}
If there exists $\eta_n\to0$ such that
\begin{equation}\label{eq:proxy-convenient}
\rho_{n,k}(\eta_n)\,\binom {d_n}{k}\to 0,\qquad
\omega_n(2\eta_n)\to0,\qquad
H_{n,k}(y)=O(1),
\end{equation}
then $D_{n,k}(y)\to0$.
\end{corollary}

\begin{proof}
In view of \eqref{eq:proxy-bound}, to prove \eqref{eq:D-bound},
it suffices to bound $A_{n,k}(\eta,y)$ and $B_{n,k}(\eta,y)$ in terms of $\rho_{n,k}(\eta)$
and $H_{n,k}(y)$:
\begin{align}
A_{n,k}(\eta,y)
&=
\sum_{|I|=k}
\Pp\{|T_n-T_{n,I}^{\circ}|>\eta\}
\,\Pp\{\cE_{n,I}(y)\} \notag\\
&\le
\rho_{n,k}(\eta)
\sum_{|I|=k}\Pp\{\cE_{n,I}(y)\}
=
\rho_{n,k}(\eta)\, H_{n,k}(y),
\label{eq:convenient-A}
\end{align}
and
\begin{align}
B_{n,k}(\eta,y)
&=
\sum_{|I|=k}
\Pp\{|T_n-T_{n,I}^{\circ}|>\eta,\cE_{n,I}(y)\} \notag\\
&\le
\sum_{|I|=k}
\Pp\{|T_n-T_{n,I}^{\circ}|>\eta\}
\le
\rho_{n,k}(\eta)\,\binom {d_n}{k}.
\label{eq:convenient-B}
\end{align}
Notice that the upper bound in \eqref{eq:convenient-B} also bounds $A_{n,k}(\eta_n,y)$,
so \eqref{eq:proxy-convenient} implies that the right-hand side of \eqref{eq:D-bound} tends to $0$
taking $\eta=\eta_n$.
\end{proof}

Although the upper bound in \eqref{eq:convenient-B} also bounds $A_{n,k}(\eta_n,y)$,
it is typically much cruder than \eqref{eq:convenient-A}. \linebreak Therefore,
one needs either a strong enough bound for $\rho_{n,k}$ to pay the cost of factor $\binom {d_n}k$,
as in \eqref{eq:proxy-convenient} and in the proofs of Sections~\ref{subsec:lixue}
and~\ref{subsec:zhaowangfeng}, or a more delicate estimate for the weighted conditional proxy error.

\begin{corollary}[$L^2$ proxy criterion]\label{cor:L2proxy}
Assuming $T_{n,I}^{\circ}\ind \cE_{n,I}(y)$, suppose for fixed $k\ge1$ and $y\in\R$ that
\begin{align}
\sup_{|I|=k}\Ee\bigl(T_n-T_{n,I}^{\circ}\bigr)^2&\to0,
\label{eq:L2ord}\\
\sum_{|I|=k}\Ee\left[
\bigl(T_n-T_{n,I}^{\circ}\bigr)^2\one\{\cE_{n,I}(y)\}
\right]&\to0,
\label{eq:L2rare}\\
H_{n,k}(y)&=O(1).
\label{eq:HO1}
\end{align}
The first condition controls the weighted unconditional proxy error,
while the second controls the weighted conditional proxy error.
If $T_n\Rightarrow T$ with a continuous distribution, then $D_{n,k}(y)\to0$.
\end{corollary}
\clearpage

\begin{proof}
Let $a_n^{\text{ord}}$ and $a_n^{\text{rare}}$ denote the two quantities in \eqref{eq:L2ord}
and \eqref{eq:L2rare} and let $a_n$ be the maximum of them.
Set $\eta_n=a_n^{1/4}$ when $a_n>0$ and $\eta_n=n^{-1}$ when $a_n=0$.
Markov's inequality gives
\vspace{-2.5pt}\[\vspace{-2.5pt}
A_{n,k}(\eta_n,y)
\le
H_{n,k}(y)\frac{a_n^{\text{ord}}}{\eta_n^2}
=O(a_n^{1/2})
\]
and
\vspace{-2.5pt}\[\vspace{-2.5pt}
B_{n,k}(\eta_n,y)
\le
\frac{a_n^{\text{rare}}}{\eta_n^2}
\le a_n^{1/2}.
\]
If $T_n\Rightarrow T$ and $T$ has a continuous distribution function $F$,
then Pólya's theorem yields
\[
\|F_n-F\|_\infty\longrightarrow 0.
\]
Since every continuous distribution function is uniformly continuous,
for every deterministic $\eta_n\to0$,
\begin{equation}\label{eq:anti}
\omega_n(2\eta_n)
\le
2\|F_n-F\|_\infty
+
\sup_{x\in\R}\{F(x+2\eta_n)-F(x-2\eta_n)\}
\longrightarrow0.
\end{equation}
Thus Theorem~\ref{thm:proxy} applies.
\end{proof}

\begin{proposition}[Multivariate $\bm T_n$ version]\label{rem:multivariate}
For fixed $K$, Theorem~\ref{thm:proxy} extends verbatim to $\bm T_n\in\R^K$. For $\eta>0$, define
\begin{align*}
A_{n,k}(\eta,y)
&=
\sum_{|I|=k}
\Pp\{\|\bm T_n-\bm T_{n,I}^{\circ}\|_\infty>\eta\}
\,\Pp\{\cE_{n,I}(y)\},
\\
B_{n,k}(\eta,y)
&=
\sum_{|I|=k}
\Pp\{\|\bm T_n-\bm T_{n,I}^{\circ}\|_\infty>\eta,\cE_{n,I}(y)\},
\\
\omega_n(\eta)
&=
\sup_{\bm x\in\R^K}
\{F_n(\bm x+\eta\bm 1)-F_n(\bm x-\eta\bm 1)\},
\end{align*}
where \(F_n(\bm x):=\Pp(\bm T_n\le\bm x)\).
Then, for every $\eta>0$,
\begin{align*}
D_{n,k}(y)
&\ =\ 
\sup_{\bm x\in \R^K}
\sum_{|I|=k}
\left|
\Pp\{\bm T_n\le\bm x,\cE_{n,I}(y)\}
-
\Pp\{\bm T_n\le\bm x\}\Pp\{\cE_{n,I}(y)\}
\right| \\
&\ \le\ 
A_{n,k}(\eta,y)
+
B_{n,k}(\eta,y)
+
\omega_n(2\eta)\times H_{n,k}(y).
\end{align*}
Theorem~\ref{thm:bonferroni} and Corollaries~\ref{cor:uniformproxy} and~\ref{cor:L2proxy} also
extend verbatim. Their proofs are unchanged. The multivariate version of
Corollary~\ref{cor:uniformproxy} is used in Section~\ref{subsec:zhaowangfeng}.
\end{proposition}

\enlargethispage*{.5\baselineskip}

\section{Local deletion}\label{sec:deletion}

The independent-proxy formulation contains local deletion as an immediate special case. Let
\[
\{X_{n,v}:v\in\cV_n\}
\]
be a triangular array with a dependency graph $G_n=(\cV_n,E_n)$.
Recall that $G_n$ is a dependency graph if whenever two disjoint vertex sets $A,B\subseteq\cV_n$
have no edge between them, the families $\{X_{n,v}:v\in A\}$ and $\{X_{n,v}:v\in B\}$ are
independent. Dependency graphs are classical tools for normal approximation of locally dependent
sums; see \citet*{Janson1988,BaldiRinott1989}.
Consider
\[
T_n
=
\frac{1}{\sigma_n}
\sum_{v\in\cV_n}
\{X_{n,v}-\Ee X_{n,v}\},
\qquad
\text{where \(\sigma_n>0\) is a normalizing factor.}
\]

\noindent
For $I\subseteq\cV_n$, let $\Gamma_n[I]$ denote the closed graph neighborhood of $I$, i.e.,
the vertices in $I$ together with every vertex adjacent to at least one vertex of $I$.
Suppose $\cE_{n,v}(y)$ is measurable with respect to $X_{n,v}$
and set $\cE_{n,I}(y)=\bigcap_{v\in I}\cE_{n,v}(y)$. Define the proxy statistic
\begin{equation}\label{eq:deleted-proxy}
T_{n,I}^{\circ}
=
\frac{1}{\sigma_n}
\sum_{v\notin\Gamma_n[I]}
\{X_{n,v}-\Ee X_{n,v}\}
\end{equation}
by the ``local deletion''. There is no edge from $I$ to $\cV_n\setminus\Gamma_n[I]$,
hence \(T_{n,I}^{\circ}\ind \cE_{n,I}(y)\). Thus Theorem~\ref{thm:proxy}
gives the following.
\begin{corollary}[Local-deletion construction]\label{cor:deletion}
Fix $k$ and $y$. Suppose that $T_{n,I}^{\circ}\ind \cE_{n,I}(y)$ and define
\begin{equation}\label{eq:deleted-rem}
R_{n,I}
:=
T_n-T_{n,I}^{\circ}
=
\frac{1}{\sigma_n}
\sum_{v\in\Gamma_n[I]}
\{X_{n,v}-\Ee X_{n,v}\}.
\end{equation}
If $H_{n,k}(y)=O(1)$, $T_n\Rightarrow T$ with continuous distribution,
and there exists $\eta_n\to0$ such that
\begin{align}
\sum_{|I|=k}\Pp(|R_{n,I}|>\eta_n)\,\Pp\{\cE_{n,I}(y)\}&\to0,
\label{eq:delordinary}\\
\sum_{|I|=k}\Pp\{|R_{n,I}|>\eta_n,\cE_{n,I}(y)\}&\to0,
\label{eq:delrare}
\end{align}
then $D_{n,k}(y)\to0$. If these conditions hold for every fixed $k$
and $\lim_{r\to\infty}\limsup_{n\to\infty}H_{n,2r}(y)=0$ holds, then $T_n$
and the max event $\cU_n(y)$ are asymptotically independent.
\end{corollary}

\begin{remark}[What a small dependency degree does and does not imply]\label{rem:degree}
Let $N_n=|\cV_n|$ and let $D_n$ be the maximum graph degree. For fixed $k$,
\( |\Gamma_n[I]|\le k(D_n+1). \) Thus $D_n=o(N_n)$ says that finitely many extreme coordinates
contaminate only a vanishing fraction of the summands. This is a useful structural heuristic,
but it is not by itself sufficient: the deleted neighborhood may carry a non-negligible fraction of
the total variance or may become large precisely on the extreme event.
Conditions \eqref{eq:delordinary}--\eqref{eq:delrare} express the missing ``local energy''
requirement.
\end{remark}

\enlargethispage*{.75\baselineskip}

\section{Examples and applications}
\label{sec:applications}

\subsection{The i.i.d. Gaussian sum and maximum}
\label{subsec:iidGaussian}

Let $Z,Z_1,\ldots,Z_n$ be i.i.d. $N(0,1)$ and define
\[
T_n=\frac{1}{\sqrt n}\sum_{i=1}^n Z_i,
\qquad
M_n=\max_{1\le i\le n}Z_i.
\]
Set $\cE_{n,i}(y)=\{Z_i>u_n(y)\}$ where $u_n(y)$ is chosen to be the Gaussian extreme-value
threshold
\[
u_n(y)
=\sqrt{2\log n}
-\frac{\log\log n+\log(4\pi)}
{2\sqrt{2\log n}}
+\frac{y}{\sqrt{2\log n}}
\]
satisfying
\[
q_n(y):=\Pp(Z>u_n(y))
=\overline\Phi(u_n(y))\sim\frac{e^{-y}}{n},\qquad
n\overline\Phi(u_n(y))\to e^{-y},
\]
with \(\overline\Phi\) the standard normal survival function.

For $I\subseteq\llbracket n\rrbracket$, define the proxy statistic and its remainder
\[
T_{n,I}^{\circ}
=
\frac{1}{\sqrt n}\sum_{i\notin I}Z_i,
\qquad
R_{n,I}
:=T_n-T_{n,I}^{\circ}
=
\frac{1}{\sqrt n}\sum_{i\in I}Z_i.
\]
This is a special case of the local-deletion mechanism of Section~\ref{sec:deletion}
when the graph is edgeless and $\Gamma_n[I]=I$. Then $T_{n,I}^{\circ}\ind \cE_{n,I}(y)$.

Fix $k\ge1$ and $y\in\R$. In this example, the proxy criterion can be verified particularly simply
through conditions \eqref{eq:L2ord}--\eqref{eq:HO1} of Corollary~\ref{cor:L2proxy}.

We first verify \eqref{eq:L2ord}. Since for any $I\subseteq\llbracket n\rrbracket$ with $|I|=k$,
\[
\Ee R_{n,I}^2=\frac{k}{n},
\]
we have
\begin{equation}\label{eq:iid-L2ord}
\sup_{|I|=k}\Ee R_{n,I}^2
=\frac{k}{n}\longrightarrow0.
\end{equation}

Now we verify \eqref{eq:L2rare}. Using the Cauchy inequality $(\sum_{i\in I}Z_i)^2\le k\sum_{i\in I}Z_i^2$
and independence,
\begin{align*}
\Ee\left[R_{n,I}^2\one\{\cE_{n,I}(y)\}\right]
&\le
\frac{k}{n}\,
\sum_{i\in I}\Ee\,\bigg[Z_i^2\one\{Z_i>u_n(y)\}\prod_{j\in
I\backslash\{i\}}\one\{Z_j>u_n(y)\}\bigg]\\
&=
\frac{k^2}{n}\,
\Ee\left[Z^2\one\{Z>u_n(y)\}\right]
q_n(y)^{k-1},
\end{align*}
where $Z\sim N(0,1)$. Hence
\begin{align}
&\sum_{|I|=k}
\Ee\left[R_{n,I}^2\one\{\cE_{n,I}(y)\}\right]
\le
\binom nk\frac{k^2}{n}\,
\Ee\left[Z^2\one\{Z>u_n(y)\}\right]
q_n(y)^{k-1}.
\label{eq:iid-L2rare}
\end{align}
For a standard Gaussian variable $Z$,
\[
\Ee\left[Z^2\one\{Z>u\}\right]
=u\phi(u)+\overline\Phi(u)
\sim u^2\overline\Phi(u),
\qquad u\to\infty.
\]
Since $u_n(y)^2\sim 2\log n$ and $nq_n(y)\to e^{-y}$,
\[
\Ee\left[Z^2\one\{Z>u_n(y)\}\right]
q_n(y)^{k-1}
\sim
u_n(y)^2 q_n(y)^{k}.
\]
Thus the right-hand side of \eqref{eq:iid-L2rare} is on the order of
\begin{equation*}
\binom nk\frac{k^2}{n}\big(\log n\big)\left(\frac{e^{-y}}{n}\right)^k \le \frac{\log
n}{n}\frac{k^2e^{-yk}}{k!}
\end{equation*}
and therefore
\begin{align}
&\sum_{|I|=k}
\Ee\left[R_{n,I}^2\one\{\cE_{n,I}(y)\}\right]
\longrightarrow 0.
\label{eq:iid-L2rare-0}
\end{align}

Last, \eqref{eq:HO1} holds because for fixed $k$,
\begin{equation}\label{eq:iid-HO1}
H_{n,k}(y)
=
\binom nk q_n(y)^k
\longrightarrow
\frac{e^{-ky}}{k!}.
\end{equation}
\clearpage

Combining \eqref{eq:iid-L2ord}, \eqref{eq:iid-L2rare-0}, and \eqref{eq:iid-HO1},
Corollary~\ref{cor:L2proxy} gives
\[
D_{n,k}(y)\longrightarrow0
\qquad\text{for every fixed }k,
\]
that is, the finite-order condition \eqref{eq:D-to-zero} in Theorem~\ref{thm:bonferroni} holds.
Since also for each fixed $r$,
\[
H_{n,2r}(y)\longrightarrow \frac{e^{-2ry}}{(2r)!}
\]
the high-order condition \eqref{eq:H-tail} in Theorem~\ref{thm:bonferroni} follows by letting
$r\to\infty$. The theorem therefore recovers the classical asymptotic independence between the
normalized Gaussian sum $T_n$ and maximum $M_n$. Corollary~\ref{cor:poisson} further gives for every
$x,y\in\R$,
\begin{equation*}
\Pp\{T_n\le x, M_n\le u_n(y)\}\longrightarrow \Phi(x)\exp\{-e^{-y}\}.
\end{equation*}

This example isolates the mechanism in its simplest form: delete the finitely many coordinates
associated with an extreme configuration, obtain exact independence for the deleted proxy,
and use an $L^2$-negligibility calculation to transfer that independence back to the original sum
statistic.

\enlargethispage*{.5\baselineskip}

\subsection{Dependent Gaussian sum and maximum}\label{subsec:dpdGau}

Let $Z_{n,1},\ldots,Z_{n,n}$ again be jointly Gaussian with $N(0,1)$ marginals.
Unlike in Section~\ref{subsec:iidGaussian}, however, we now allow them to be \emph{dependent}.
The dependency-graph specialization introduced in Section~\ref{sec:deletion}
and used in Section~\ref{subsec:iidGaussian} is no longer suitable,
since nonzero correlations may occur between all pairs of variables.
Nevertheless, Theorem~\ref{thm:proxy} remains applicable when the proxy is constructed via Gaussian
regression and residualization rather than local deletion.

Let $\bm Z_n=(Z_{n,1},\ldots,Z_{n,n})^\top\sim N(0,\mathbf\Sigma_n)$,
where $\mathbf\Sigma_n=(\sigma_{ij})_{1\le i,j\le n}$.
Denoting $S_n=\sum_{i=1}^n Z_{n,i}$ and $\sigma_n^2=\operatorname{Var}(S_n)$, we consider
the sum and maximum
\[
T_n=S_n/\sigma_n,
\qquad
M_n=\max_{1\le i\le n}Z_{n,i}.
\]
Set once again the classical Gaussian extreme-value threshold
\[
u_n(y)
=\sqrt{2\log n}
-\frac{\log\log n+\log(4\pi)}
{2\sqrt{2\log n}}
+\frac{y}{\sqrt{2\log n}}
\]
and $\cE_{n,i}(y)=\{Z_{n,i}>u_n(y)\}$. We show that, for every $x,y\in\R$,
\begin{equation}\label{eq:mq-joint}
\Pp\{T_n\le x,\ M_n\le u_n(y)\}
\longrightarrow\Phi(x)\exp\{-e^{-y}\};
\end{equation}
that is, $T_n$ and $M_n$ are asymptotically independent,
with exactly the same joint limit 
as in Section~\ref{subsec:iidGaussian}. 

To facilitate our analysis, we need to impose certain assumptions on the data generation process.
Suppose that
\begin{equation*}
\mathbf\Sigma_n\succeq0,
\qquad
\sigma_{ii}=1,
\qquad
\bm 1_n^\top\mathbf\Sigma_n\bm 1_n>0,
\end{equation*}
which is assumed only for standardization and, in particular, $\bm 1_n^\top\mathbf\Sigma_n\bm 1_n>0$
is needed only to define $T_n$; also suppose that, for fixed $\varrho\in(0,1)$ and $\xi>0$,
\begin{equation}\label{eq:mq-assumptions}
\max_{i\ne j}|\sigma_{ij}|\le\varrho,
\qquad
\max_i\sum_{j=1}^n\sigma_{ij}^2\le(\log n)^\xi.
\end{equation}
No eigenvalue or condition-number assumption is imposed.

\newpage
\paragraph{Nearly uniform proxy criterion.}
We follow the approach of Corollary~\ref{cor:uniformproxy}. 
Under the present assumptions, however, establishing a uniform bound of the form
\begin{equation*}
\rho_{n,k}(\eta)
:=\sup_{|I|=k}\Pp(|T_n-T_{n,I}^{\circ}|>\eta)
\end{equation*}
may be difficult.
We therefore partition the sets $I$ into two classes, termed \emph{good} and \emph{bad}.
Their precise definitions will be given below; for now, the key point is that the good sets are
well behaved and make the dominant contribution.

\enlargethispage*{.5\baselineskip}

We write
\begin{align*}
H_{n,k}(y)
&=
\sum_{|I|=k}
\Pp\{\cE_{n,I}(y)\}\\
&=
\sum_{\substack{|I|=k\\I\text{ good}}}
\Pp\{\cE_{n,I}(y)\}
+
\sum_{\substack{|I|=k\\I\text{ bad}}}
\Pp\{\cE_{n,I}(y)\}
=:H_{n,k}^{\mathrm{good}}(y)+H_{n,k}^{\mathrm{bad}}(y),
\end{align*}
and we define
\begin{equation*}
\rho_{n,k}^{\mathrm{good}}(\eta_n)
:=\sup_{\substack{|I|=k,\ I\text{ good}}}\Pp(|T_n-T_{n,I}^{\circ}|>\eta_n).
\end{equation*}
Using the definitions of the two proxy errors $A_{n,k}(\eta_n,y)$ and $B_{n,k}(\eta_n,y)$ in
Theorem~\ref{thm:proxy}, we have
\begin{align}
A_{n,k}(\eta_n,y)
&=
\sum_{|I|=k}
\Pp\{|T_n-T_{n,I}^{\circ}|>\eta_n\}
\,\Pp\{\cE_{n,I}(y)\} \notag\\
&=
\sum_{\substack{|I|=k\\I\text{ good}}}
\Pp\{|T_n-T_{n,I}^{\circ}|>\eta_n\}
\,\Pp\{\cE_{n,I}(y)\}
+
\sum_{\substack{|I|=k\\I\text{ bad}}}
\Pp\{|T_n-T_{n,I}^{\circ}|>\eta_n\}
\,\Pp\{\cE_{n,I}(y)\} \notag\\
&\le
\rho_{n,k}^{\mathrm{good}}(\eta_n)
\sum_{\substack{|I|=k\\I\text{ good}}}
\Pp\{\cE_{n,I}(y)\}
+
\sum_{\substack{|I|=k\\I\text{ bad}}}
\Pp\{\cE_{n,I}(y)\} \notag\\
&\le
\rho_{n,k}^{\mathrm{good}}(\eta_n)\, H_{n,k}^{\mathrm{good}}(y)
+
H_{n,k}^{\mathrm{bad}}(y),
\notag
\end{align}
and
\begin{align}
B_{n,k}(\eta_n,y)
&=
\sum_{|I|=k}
\Pp\{|T_n-T_{n,I}^{\circ}|>\eta_n,\cE_{n,I}(y)\} \notag\\
&=
\sum_{\substack{|I|=k\\I\text{ good}}}
\Pp\{|T_n-T_{n,I}^{\circ}|>\eta_n,\cE_{n,I}(y)\}
+
\sum_{\substack{|I|=k\\I\text{ bad}}}
\Pp\{|T_n-T_{n,I}^{\circ}|>\eta_n,\cE_{n,I}(y)\} \notag\\
&\le
\sum_{\substack{|I|=k\\I\text{ good}}}
\Pp\{|T_n-T_{n,I}^{\circ}|>\eta_n\}
+
\sum_{\substack{|I|=k\\I\text{ bad}}}
\Pp\{\cE_{n,I}(y)\}
\notag\\
&\le
\rho_{n,k}^{\mathrm{good}}(\eta_n)\,\binom {d_n}{k}
+
H_{n,k}^{\mathrm{bad}}(y).
\notag
\end{align}
Therefore, if one can find proxy $T_{n,I}^{\circ}\ind \cE_{n,I}(y)$
and prove for some $\eta_n\to0$ that
\begin{equation}\label{eq:dpdGau-proxy-convenient}
\rho_{n,k}^{\mathrm{good}}(\eta_n)\,\binom {d_n}{k}\to 0,\qquad
\omega_n(2\eta_n)\to0,\qquad
H_{n,k}^{\mathrm{good}}(y)=O(1),\qquad
H_{n,k}^{\mathrm{bad}}(y)\to 0,
\end{equation}
then Theorem~\ref{thm:proxy} gives $D_{n,k}(y)\to0$.

\newpage
\paragraph{Gaussian regression and residualization.}
Fix $k\ge1$ and $y\in\R$. For every $I$ with $|I|=k$, define the sum proxy and its remainder by
\begin{equation}\label{eq:dpdGau-decomp1}
T_{n,I}^{\circ}:=T_n-\Ee(T_n\mid\bm Z_{n,I}),
\qquad
R_{n,I}:=T_n-T_{n,I}^{\circ}.
\end{equation}
This decomposition remains valid even when $\mathbf\Sigma_{n,I}$ is singular,
since the Gaussian conditional expectation is well defined. Let \(I=\{i_1,\ldots,i_k\}\) and denote
\[
\bm Z_{n,I}=(Z_{n,i_1},\ldots,Z_{n,i_k})^\top,\qquad
\mathbf\Sigma_{n,I}=\operatorname{Var}(\bm Z_{n,I})=(\mathbf\Sigma_{n})_{I,I},
\]
and define the regression coefficient vector
\[
\bm r_{n,I}
:=\operatorname{Cov}(\bm Z_{n,I},T_n)
=\frac{(\mathbf\Sigma_n\bm1_n)_I}{\sigma_n}.
\]
The decomposition in \eqref{eq:dpdGau-decomp1} can be explicitly written as
\begin{equation}\label{eq:dpdGau-decomp2}
T_{n,I}^{\circ}=T_n-\bm r_{n,I}^\top\mathbf\Sigma_{n,I}^{-1}\bm Z_{n,I},
\qquad
R_{n,I}
=\bm r_{n,I}^\top\mathbf\Sigma_{n,I}^{-1}\bm Z_{n,I}
\end{equation}
when $\mathbf\Sigma_{n,I}$ is nonsingular, with the inverse $\mathbf\Sigma_{n,I}^{-1}$ replaced by
Moore–Penrose pseudoinverse when $\mathbf\Sigma_{n,I}$ is singular;
the proxy $T_{n,I}^{\circ}$ is indeed the Gaussian residualization. Joint Gaussianity gives
\[
T_{n,I}^{\circ}\ind\bm Z_{n,I},
\qquad\text{hence}\quad
T_{n,I}^{\circ}\ind\cE_{n,I}(y).
\]

\bigskip
\bigskip

\paragraph{Good and bad sets $I\subseteq\llbracket n\rrbracket$.}
We now bound the quantities associated with the good and bad $I$ required by \eqref{eq:dpdGau-proxy-convenient}
for the proxy and remainder in the form of \eqref{eq:dpdGau-decomp2}, starting with the definitions of
good and bad $I$. \looseness=-1

Fix $k\ge1$ and $y\in\R$. Consider every $I\subseteq\llbracket n\rrbracket$ with $|I|=k$.
Take $$\delta_n=(\log n)^{-2}.$$ We call a pair $i\ne j$ \emph{weakly correlated}
if $|\sigma_{ij}|\le\delta_n$ and \emph{strongly correlated} otherwise.
We call a set~$I$ \emph{weakly correlated} if $|\sigma_{ij}|\le\delta_n$ for all distinct
$i,j\in I$.

We call a set $I$ \emph{good} if it satisfies \emph{both} of the following conditions:
\begin{itemize}[itemsep=-.1em, topsep=.5em]
\item \emph{$I$ itself is weakly correlated}: $|\sigma_{ij}|\le\delta_n$ for all distinct
$i,j\in I$, and
\item \emph{$\bm r_{n,I}$ is uniformly weak}: $|r_{n,i}|\le\delta_n$ for all $i\in I$,
or alternatively, \(\lVert\bm r_{n,I}\rVert_{\infty} \le \delta_n\).
\end{itemize}
Every \emph{good} set $I$ is \emph{weakly correlated}, but the converse cannot be guaranteed.
\emph{Bad} sets $I$ are of two types:
\begin{itemize}[itemsep=-.1em, topsep=.5em]
\item \emph{Type I bad (strong)}: there exists at least one pair of distinct $i,j\in I$ such that
\(|\sigma_{ij}|>\delta_n\).
\item \emph{Type II bad (weak)}: $I$ is weakly correlated, but
\(\lVert\bm r_{n,I}\rVert_{\infty} > \delta_n\), that is, there exists $i\in I$ such that
\[
|r_{n,i}|=|(\bm r_n)_i|
=\left|\frac{\sum_{j=1}^n\sigma_{ij}}{\sigma_n}\right|
=\left|\frac{\sum_{j=1}^n\sigma_{ij}}{\sqrt{\sum_{i,j}\sigma_{ij}}}\right|
>\delta_n,
\]
where
\[
\bm r_n:=\operatorname{Cov}(\bm Z_n,T_n)
=\frac{\mathbf\Sigma_n\bm 1_n}{\sigma_n}.
\]
\end{itemize}

\newpage
\paragraph{Compute $\Pp\{\cE_{n,I}(y)\}$ for all weakly correlated $I$.}
$I$ is weakly correlated if $|\sigma_{ij}|\le\delta_n$ for all distinct $i,j\in I$.
By the Gershgorin circle theorem \citep*{Gerschgorin1931}, the following squeeze argument (squeeze
principle) shows that $\Pp\{\cE_{n,I}(y)\}$ behaves asymptotically as if $\Sigma_{n,I}=\mathbf I_k$.

Setting $a_{k,n}=(k-1)\delta_n$, the eigenvalues of $\mathbf\Sigma_{n,I}$ lie between $1-a_{k,n}$
and $1+a_{k,n}$ and
\[
(1-a_{k,n})\,\mathbf{I}_k\preceq \mathbf\Sigma_{n,I}\preceq(1+a_{k,n})\,\mathbf{I}_k.
\]
Consequently, for all sufficiently large $n$,
\[
\frac{\mathbf{I}_k}{1+a_{k,n}}\preceq \mathbf\Sigma_{n,I}^{-1}\preceq\frac{\mathbf{I}_k}{1-a_{k,n}},
\qquad
(1-a_{k,n})^k\le\det \mathbf\Sigma_{n,I}\le(1+a_{k,n})^k.
\]
Applying these bounds directly to the Gaussian density for $\bm Z_{n,I}$
and integrating over the region \( (z,+\infty)^k \) for $z>0$, we obtain
\[
\begin{aligned}
\left(\frac{1-a_{k,n}}{1+a_{k,n}}\right)^{k/2}
\left[\overline\Phi\left(\frac{z}{\sqrt{1-a_{k,n}}}\right)\right]^k
&\le\ 
\Pp( Z_{n,i_1} >z,\ldots, Z_{n,i_k} >z)\\
&\hspace{-5em}=\ 
\frac{1}{(2\pi)^{k/2}\det(\mathbf\Sigma_{n,I})^{1/2}}
\int_{z}^{+\infty}\cdots\int_{z}^{+\infty}
\exp\left(-\frac12\bm x_{I}^{\top}\mathbf\Sigma_{n,I}^{-1}\bm x_{I}\right)
\,{\rm d}x_{i_1}\cdots{\rm d}x_{i_k}\\
&\hspace{-5em}\le\ 
\left(\frac{1+a_{k,n}}{1-a_{k,n}}\right)^{k/2}
\left[\overline\Phi\left(\frac{z}{\sqrt{1+a_{k,n}}}\right)\right]^k.
\end{aligned}
\]
Taking \( z = u_n(y) \), since
\[
u_n(y)^2
=2\log n-\log\log n-\log(4\pi)+2y
+O\left(\frac{(\log\log n)^2}{\log n}\right)
\]
and $a_{k,n} u_n(y)^2=O(k\delta_n \log n)=O\{k(\log n)^{-1}\}=o(1)$ for any fixed $k$,
Mills' ratio yields
\begin{align*}
\frac{\overline\Phi(u_n(y)/\sqrt{1\pm a_{k,n}})}
{\overline\Phi(u_n(y))}
&\sim
\frac{\phi(u_n(y)/\sqrt{1\pm a_{k,n}})}
{\phi(u_n(y))}
\times
\frac{u_n(y)}
{u_n(y)/\sqrt{1\pm a_{k,n}}}\\
&=
\sqrt{1\pm a_{k,n}}
\exp\left\{\frac{z^2}{2}
\left(1-\frac1{1\pm a_{k,n}}\right)\right\}
\longrightarrow1.
\end{align*}
Meanwhile,
\[
q_n(y):=\overline\Phi(u_n(y))\sim \frac{e^{-y}}{n}.
\]
Therefore, the preceding bounds imply
\begin{equation}\label{eq:mq-weak-tail}
\Pp\{\cE_{n,I}(y)\}
=\Pp( Z_{n,i_1} >u_n(y),\ldots, Z_{n,i_k} >u_n(y))
=(1+o(1))q_n(y)^k
=(1+o(1))\left(\frac{e^{-y}}n\right)^k.
\end{equation}
All error bounds depend only on the fixed $k$ and $\delta_n=(\log n)^{-2}$,
and not on the choice of index set $I=\{i_1,\ldots,i_k\}$.
Hence, the convergence is uniformly over all \emph{good} (in fact, over all
\emph{weakly~correlated})~$I$, completing the proof.

\newpage
\paragraph{Compute $H_{n,k}^{\mathrm{good}}(y)$.}
Recall the definition
\[
H_{n,k}^{\mathrm{good}}(y) := \sum_{\substack{|I|=k:\ I\text{ good}}}\Pp\{\cE_{n,I}(y)\}
\]
We will show below that the cardinality $\#\{|I|=k:\ I\text{ bad}\}$ is $o(n^k)$:
\begin{equation}\label{eq:dpdGau-counts-bad}
\sum_{\substack{|I|=k:\ I\text{ bad}}}1 = o(n^k)
\end{equation}
Then \eqref{eq:mq-weak-tail} gives
\begin{equation}\label{eq:dpdGau-H-good}
H_{n,k}^{\mathrm{good}}(y) = \sum_{\substack{|I|=k:\ I\text{ good}}}\Pp\{\cE_{n,I}(y)\} =
(1+o(1))\left(\frac{e^{-y}}n\right)^k\left\{\binom nk - o(n^k)\right\} = O(1).
\end{equation}

\bigskip
\bigskip

\enlargethispage*{.5\baselineskip}

\paragraph{Bound $\rho_{n,k}^{\mathrm{good}}(\eta_n)$.}
Recall the definition of $\rho_{n,k}^{\mathrm{good}}(\eta_n)$
\begin{equation*}
\rho_{n,k}^{\mathrm{good}}(\eta_n)
=\sup_{\substack{|I|=k\\I\text{ good}}}\Pp(|T_n-T_{n,I}^{\circ}|>\eta_n)
=\sup_{\substack{|I|=k\\I\text{ good}}}\Pp(|R_{n,I}|>\eta_n)
\end{equation*}
and the explicit form of our constructed remainder $R_{n,I}$
\[
R_{n,I}
=\bm r_{n,I}^{\top} \mathbf\Sigma_{n,I}^{-1}\bm Z_{n,I},
\qquad\text{with}\quad
\mathbb E R_{n,I}^2
=\bm r_{n,I}^{\top}\mathbf\Sigma_{n,I}^{-1}\bm r_{n,I}.
\]
Notice that the covariance matrix $\mathbf\Sigma_{n,I}$ is always invertible for good (and thus
weakly correlated) $I$, since its eigenvalues all lie between $1-a_{k,n}$ and $1+a_{k,n}$.
Since \(R_{n,I}\) is a linear transformation of a Gaussian vector
\[
R_{n,I}\sim
N\left(
0,\,
\bm r_{n,I}^\top
\mathbf\Sigma_{n,I}^{-1}
\bm r_{n,I}
\right).
\]
For good \(I\), we have \(\lambda_{\min}(\mathbf\Sigma_{n,I}) \ge1-(k-1)\delta_n\) by its weak
correlation within $I$, and
\[
\|\bm r_{n,I}\|^2\le k\delta_n^2,
\]
by its weak correlation with the standardized sum. Thus
\[
\mathbb E R_{n,I}^2
\le
\frac{k\delta_n^2}{1-(k-1)\delta_n}.
\]
The Gaussian tail bound (Chernoff bound) gives
\[
\Pp(|R_{n,I}|>\eta_n)
\le
2\exp\left\{
-\frac{\eta_n^2[1-(k-1)\delta_n]}
{2k\delta_n^2}
\right\}.
\]
Taking $\eta_n=(\log n)^{-1/2}$ and recalling $\delta_n=(\log n)^{-2}$ yields
\[
\sup_{\substack{|I|=k\\ I\text{ good}}}
\Pp\big\{|R_{n,I}|>(\log n)^{-1/2}\big\}
\le
2\exp\left\{
-\frac{(1+o(1))(\log n)^3}{2k}
\right\}.
\]
This convergence rate is fast enough to offset the combinatorial factor $\binom {d_n}{k}$.

\newpage
\paragraph{Bound $H_{n,k}^{\mathrm{bad}}(y)$.}
Recalling the two types of bad sets $I$:
\begin{itemize}[itemsep=-.1em, topsep=.5em]
\item \emph{Type I bad (strong)}: there exists at least one pair of distinct $i,j\in I$ such that
\(|\sigma_{ij}|>\delta_n\)
\item \emph{Type II bad (weak)}: I is weakly correlated, however, there exists $i\in I$ such that \[
|r_{n,i}|=\left|\frac{\sum_{j=1}^n\sigma_{ij}}{\sigma_n}\right|
>\delta_n.
\]
\end{itemize}
We further write
\[
H_{n,k}^{\mathrm{bad}}(y) = H_{n,k}^{\mathrm{bad~(strong)}}(y) + H_{n,k}^{\mathrm{bad~(weak)}}(y)
:=\sum_{\substack{|I|=k\\I\text{ bad (strong)}}}
\Pp\{\cE_{n,I}(y)\}
+
\sum_{\substack{|I|=k\\I\text{ bad (weak)}}}
\Pp\{\cE_{n,I}(y)\}.
\]

\bigskip
\bigskip

\enlargethispage*{.5\baselineskip}

\noindent\underline{\underline{\itshape Control $H_{n,k}^{\mathrm{bad~(strong)}}(y)$.}}\hspace{.75em}
When $k=2$, to control
$H_{n,k}^{\mathrm{bad~(strong)}}(y)$, 
applying the normal comparison inequality in \citet*{LiShao2002} yields the desired bound.
However, this argument does not directly extend to $k\ge 3$.
We instead consider the following connected-component argument, which applies to all $k\ge 2$.
The heuristic is to partition the variables into components
such that variables within each component may be strongly correlated,
whereas dependence between different components is weak. Consequently,
the component averages are nearly independent, so their weighted sum behaves approximately like a
sum of independent random variables. 
This approach resembles classical dependency-graph arguments, 
except that distinct components here are only weakly correlated rather than independent. 

\medskip

\uwave{Strong-correlation graph and connected components.}\hspace{.55em}
We consider a graph $G_n$ with vertex set $\llbracket n\rrbracket$, \linebreak where two distinct vertices $i$
and $j$ are connected whenever $|\sigma_{ij}|>\delta_n$---they have ``strong correlation''. This graph has maximum degree at most
\begin{equation}\label{eq:dpdGau-D_n}
\#\{j\ne i:|\sigma_{ij}|>\delta_n\}\le\frac{\max_i\sum_{j=1}^n\sigma_{ij}^2}{\delta_n^2}\le(\log
n)^{\xi+4}=:D_n,
\end{equation}
where $\max_i\sum_{j=1}^n\sigma_{ij}^2\le(\log n)^\xi$ is in our assumption
\eqref{eq:mq-assumptions}.

For every Type I \emph{bad} set \(I\) with $|I|=k$, we must have \(k\ge2\)
since there exists at least one pair of distinct $i,j\in I$ such that  \(|\sigma_{ij}|>\delta_n\).
Its induced subgraph \(G_n[I]\) contains at least one edge, because such \emph{strongly correlated}
pair exists.

Let \(m\in\llbracket k-1\rrbracket\) denote the number of connected components of \(G_n[I]\)---here $m$ is
at most $k-1$ since \(G_n[I]\) contains at least one edge.
We show the number of such sets with a given $m$ is at most $C_k n^m D_n^{k-m}$ for some constant
$C_k$ depending only on $k$: note that every such \(G_n[I]\) contains a~spanning forest consisting
of \(m\) trees with each connected component containing a spanning tree.
The counting is conducted by choosing one root per component and adding the remaining vertices along
spanning trees. Each \(k\)-set is counted at least \(k!\) times, so we divide the resulting count by
\(k!\):
\begin{align*}
\#\{I\subseteq\llbracket n\rrbracket: |I|=k,\ G_n[I]\text{ has }m\text{ components}\}
&\le \frac{|\mathcal F_{k,m}|}{k!}n^mD_n^{k-m},
\end{align*}
where \(\mathcal F_{k,m}\) denotes the collection of labeled forests on index set
\(\llbracket k\rrbracket\) with exactly \(m\) connected components.
Every such labeled forest has exactly $k-m$ edges. Thus, by 
\citet*[Chapter~15]{Harary1969}
\[
|\mathcal F_{k,m}|\le \binom{\binom{k}2}{k-m}\le 2^{\binom{k}{2}}.
\]

\newpage
\uwave{Gaussian tail bound via weighted component averages.}\hspace{.75em} For any such set
$I\subseteq\llbracket n\rrbracket$ with $|I|=k$ whose induced subgraph $G_n[I]$ has $m$ components,
write its components as $\cC_1,\ldots,\cC_m$, with sizes $\kappa_{\ell}=|\cC_{\ell}|$, and define
\[
Y_\ell=\frac1{\kappa_\ell}\sum_{i\in \cC_{\ell}}Z_{n,i},
\qquad
w_\ell=\frac{\kappa_\ell}{1+(\kappa_{\ell}-1)\varrho},
\qquad
\widehat Y=\sum_{\ell=1}^m w_{\ell}Y_{\ell},
\qquad w=\sum_{\ell=1}^m w_{\ell}.
\]
Here $w_\ell$ is chosen as the reciprocal of the following upper bound on $\operatorname{Var}(Y_\ell)$:
\[
\operatorname{Var}(Y_\ell)
\le\frac1{\kappa_\ell^2}
\left(\sum_{i\in\mathcal C_\ell}\sigma_{ii}+\sum_{i,j\in\mathcal C_\ell: i\ne j}|\sigma_{ij}|\right)
\le
\frac{\kappa_\ell+\kappa_\ell(\kappa_\ell-1)\varrho}{\kappa_\ell^2}
=\frac{1+(\kappa_\ell-1)\varrho}{\kappa_\ell}.
\]
Accordingly, \(w_\ell\) can be viewed as the ``effective number of independent variables''
(a.k.a. the ``degrees of freedom'') that the component $\cC_\ell$ contributes to the Gaussian tail exponent.
For example,
\begin{itemize}[itemsep=-.1em, topsep=.5em]
\item[--] If $\kappa_\ell=1$, then $w_\ell=1$; if $\kappa_\ell\ge2$, then $1<w_\ell<\kappa_\ell$, and $w_\ell$ is strictly increasing in $\kappa_\ell$ for fixed $\varrho$. \looseness=-1
\item[--] If \(\varrho\downarrow0\), then \(w_\ell\uparrow\kappa_\ell\); the variables in this component are i.i.d.
Gaussian in the limiting case.
\item[--] If \(\varrho\uparrow1\), then \(w_\ell\downarrow1\); the variables in this component can be
up to almost fully correlated.
\end{itemize}
For any Type I bad $k$-set $I$, since at least one component $\cC_{\ell^*}$ contains two or more vertices, we have
\[
w-m=\sum_{\ell=1}^m\frac{(\kappa_{\ell}-1)(1-\varrho)}{1+(\kappa_\ell-1)\varrho}
\ge\frac{(\kappa_{\ell^*}-1)(1-\varrho)}{1+(\kappa_{\ell^*}-1)\varrho}
\ge\frac{1-\varrho}{1+\varrho}=:\gamma>0,
\]
and the variance bounds $\operatorname{Var}(Y_{\ell})\le w_{\ell}^{-1}$
and $|\operatorname{Cov}(Y_{\ell},Y_{h})|\le\delta_n$ for $\ell\ne h$ imply
\[
\operatorname{Var}(\widehat Y)
=\sum_\ell w_\ell^2\operatorname{Var}(Y_\ell)
+2\sum_{\ell<h}w_\ell w_h\operatorname{Cov}(Y_\ell,Y_h)
\le w+2\delta_n\sum_{\ell<h}w_\ell w_h
\le w+\delta_nw^2.
\]
On $\cE_{n,I}(y)$ we have $Y_{\ell}>u_n(y)$ for all $\ell$ and thus $\widehat Y>wu_n(y)$;
hence the Gaussian tail bound gives
\[
\begin{aligned}
\Pp(\mathcal E_{n,I}(y))
&\le\Pp(\widehat Y>wu_n(y))
\le\exp\left\{
-\frac{w^2u_n(y)^2}{2\operatorname{Var}(\widehat Y)}
\right\}\\
&\le\exp\left\{
-\frac{wu_n(y)^2}{2(1+\delta_n w)}
\right\}
\le C_{k,y}(\log n)^{w/2}n^{-w}
\le C_{k,y}(\log n)^{k/2}n^{-m-\gamma}.
\end{aligned}
\]
where \(m\le w\le k\) and recall that
\[
u_n(y)^2
=2\log n-\log\log n-\log(4\pi)+2y
+O\left(\frac{(\log\log n)^2}{\log n}\right).
\]
Here $\delta_n w$ can be absorbed by some constant $C_{k,y}$
because $w u_n(y)^2\times \delta_n w=O(k^2\delta_n \log n)=O\{k^2(\log n)^{-1}\}=o(1)$.

\medskip

\uwave{Combining component counts and weighted-sum tail bounds.}\hspace{.75em}
Summing over $m\in\{1,\ldots,k-1\}$ shows that the total contribution is
\begin{align*}
H_{n,k}^{\mathrm{bad~(strong)}}(y)
&=\sum_{\substack{|I|=k:\ I\text{ bad (strong)}}}
\Pp\{\cE_{n,I}(y)\}\\
&\le C_{k,y}(\log n)^{k/2}\sum_{m=1}^{k-1}D_n^{k-m}n^{-\gamma}
 \le C_{k,y}(\log n)^{c_{k,\xi}}n^{-\gamma}=o(1).\yestag\label{eq:dpdGau-H-bad-Type-I}
\end{align*}

To conclude the part of controlling $H_{n,k}^{\mathrm{bad~(strong)}}(y)$,
note also that the total number of Type I bad $k$-sets is at most
\begin{equation}\label{eq:dpdGau-counts-bad-Type-I}
\sum_{m=1}^{k-1}\binom{\binom{k}2}{k-m}n^mD_n^{k-m} \le 2^{\binom{k}2}n^{k-1}D_n = o(n^k),
\end{equation}
which will be used later.

\bigskip
\bigskip

\noindent\underline{\underline{\itshape Control $H_{n,k}^{\mathrm{bad~(weak)}}(y)$.}}\hspace{.75em}
Now it remains to handle the Type II bad sets $I$. By the spectral theorem,
$\mathbf\Sigma_n^2\preceq\lambda_{\max}(\mathbf\Sigma_n)\mathbf\Sigma_n$
and together with our assumption $\max_i\sum_{j=1}^n\sigma_{ij}^2\le(\log n)^\xi$ in
\eqref{eq:mq-assumptions},
\[
\|\bm r_n\|^2\le\lambda_{\max}(\mathbf\Sigma_n)
\le\sqrt{\operatorname{tr}(\mathbf\Sigma_n^2)}
\le\sqrt n\,(\log n)^{\xi/2}.
\]
Thus the exceptional set
\[
K_n:=\{i:|r_{n,i}|>\delta_n\}
\quad\text{satisfies}\quad
|K_n|\le\sqrt n\,(\log n)^{\xi/2+4}=o(n).
\]
Then the number of Type II \emph{bad} $k$-sets, i.e., that of \(k\)-sets containing at least one
exceptional coordinate is at most
\begin{equation}\label{eq:dpdGau-counts-bad-Type-II}
|K_n|\,\binom{n-1}{k-1}
=o(n^k)
\end{equation}
Notice that Type II bad sets $I$ are still weakly correlated,
so together with \eqref{eq:mq-weak-tail},
\begin{equation}\label{eq:dpdGau-H-bad-Type-II}
H_{n,k}^{\mathrm{bad~(weak)}}(y) = \sum_{\substack{|I|=k:\ I\text{ bad (weak)}}}\Pp\{\cE_{n,I}(y)\}
= (1+o(1))\left(\frac{e^{-y}}n\right)^k o(n^k) = o(1).
\end{equation}

To conclude the part of controlling $H_{n,k}^{\mathrm{bad~(weak)}}(y)$,
combining \eqref{eq:dpdGau-counts-bad-Type-II} and \eqref{eq:dpdGau-counts-bad-Type-I} proves
\eqref{eq:dpdGau-counts-bad}.

\bigskip
\bigskip

\paragraph{Combine good and bad parts.}
We have verified all conditions in \eqref{eq:dpdGau-proxy-convenient} except that
$\omega_n(2\eta_n)\to0$, which is obvious: $T_n$ under our normalization is exactly $N(0,1)$,
which implies
\begin{equation}\label{eq:dpdGau-omega}
\omega_n(2\eta_n)\le\frac{4\eta_n}{\sqrt{2\pi}}\longrightarrow0.
\end{equation}
Therefore, Theorem~\ref{thm:proxy} gives $D_{n,k}(y)\to0$ for every fixed~$k$.
Putting \eqref{eq:dpdGau-H-good}, \eqref{eq:dpdGau-H-bad-Type-I},
\eqref{eq:dpdGau-H-bad-Type-II} together implies
\begin{equation}\label{eq:mq-H}
H_{n,k}(y)=\sum_{|I|=k}\Pp\{\cE_{n,I}(y)\}
=\left(\frac{e^{-y}}n\right)^k \binom nk + o(1)
\longrightarrow \frac{e^{-ky}}{k!}
\qquad(k\ge1\text{ fixed}).
\end{equation}
Since \eqref{eq:mq-H} also implies the high-order binomial-moment condition \eqref{eq:H-tail},
Corollary~\ref{cor:poisson} yields \eqref{eq:mq-joint}. Thus the Gaussian sum
and maximum are asymptotically independent under dependence structure characterized by
\eqref{eq:mq-assumptions}. 

\bigskip

\citet*{McCormickQi2000} studied the joint limit of such statistics using a maximum-proxy argument.
Unlike them, the proof above directly applies Theorem~\ref{thm:proxy} by constructing a proxy for
the sum $T_n$. The assumptions are first two parts of Assumption~3 in \citet*{FengEtAl2024}.
\looseness=-1

\newpage
\subsection{High-dimensional independence testing based on pairwise sample covariances}
\label{subsec:lixue}

We next explain how the results of \citet*{LiXue2015} fit into our framework.
Their setting concerns high-dimensional mutual-independence testing 
based on pairwise sample covariances or Spearman's~$\rho$.
We illustrate how our framework applies to the former; indeed, it can accommodate both.

Write the $p_n$ data columns as $\bm X_1,\ldots,\bm X_{p_n}$ and index the $d_n=\binom {p_n}2$ pairs
in $\{(j,h):1\le j<h\le p_n\}$ in any fixed order (e.g., lexicographic order) by writing
\[
\{(j,h):1\le j<h\le p_n\} = \{(j_\ell,h_\ell): 1\le\ell\le d_n\}.
\]
Their quadratic statistic is, after centering and normalization,
\[
T_{n}
=
b_{n}\sum_{\ell=1}^{d_n}(V_\ell-1),
\qquad\text{where}\quad
V_\ell
=
\frac{(\bm X_{j_\ell}^{\top}\bm X_{h_\ell})^2}{n},
\]
and $b_n$ is some regime-dependent normalization factor defined in their Lemma~1.
The max statistic is defined as $\max_{\ell}V_\ell$.
Set $\cE_{n,\ell}(y)=\{V_\ell>u_n(y)\}$ and \(\cE_{n,I}(y)=\bigcap_{\ell\in I}\cE_{n,\ell}(y)\) for the threshold
\begin{equation}\label{eq:squareduny}
u_n(y)=4\log p_n-\log\log p_n+y.
\end{equation}

Suppose that the underlying data-generating process $\{X_{ij}:1\le i\le n,\ 1\le j\le p_n\}$ are i.i.d. with a fixed distribution satisfying
\begin{equation}\label{eq:pearson-assumptions}
\Ee X_{11}=0,\qquad
\Ee X_{11}^2=1,\qquad
\Ee X_{11}^4>1,\qquad
\Ee e^{t_0 X_{11}^2}<\infty
\quad\text{for some }t_0>0.
\end{equation}
We further assume the following growth regime:
\begin{equation}\label{eq:growing-regime}
n,p_n\to\infty,\qquad p_n\ge n,\qquad
\log p_n=o(n^{1/6}).
\end{equation}

For $I\subseteq\llbracket d_n\rrbracket$ with $|I|=k$, \citeauthor*{LiXue2015} define $W_I$ to be the
set of indices $\ell$ such that the pair $(j_\ell,h_\ell)$ shares at least one coordinate with some
pair indexed by $I$. The cardinality of $W_I$ satisfies
\[
|W_I|\le 2k p_n.
\]
Because the original columns $\bm X_1,\ldots,\bm X_{p_n}$ are independent under the null
hypothesis, the variables  $\{V_\ell:\ell\in I\}$ are independent of the collection
$\{V_\ell:\ell\notin W_I\}$. Hence the proxy with its remainder \looseness=-1
\vspace{+2.5pt}\begin{equation}
\label{eq:lixue-proxy}
T_{n,I}^{\circ}
=
b_{n}\sum_{\ell\notin W_I}(V_\ell-1),
\qquad
R_{n,I}
:=T_{n}-T_{n,I}^{\circ}
=
b_{n}\sum_{\ell\in W_I}(V_\ell-1)
\end{equation}
satisfies $T_{n,I}^{\circ}\ind\cE_{n,I}(y)$.
This is precisely the local-deletion mechanism of Section~\ref{sec:deletion}.

\enlargethispage*{.5\baselineskip}

In this example, we will use Corollary~\ref{cor:uniformproxy} to verify the proxy criterion.
We only need to verify the three conditions in \eqref{eq:proxy-convenient},
starting from the one about $\rho_{n,k}(\eta_n)$. Equation (23) of \citet*{LiXue2015} provides a uniform tail bound for the proxy remainder. Set
\[
\eta_{n}=\min\{p_n^{-1/5},\sqrt{n/p_n}\}.
\]
For each fixed $k$, their argument yields constants $C_k,c_k>0$ such that
\begin{equation}\label{eq:lixue-rho}
\rho_{n,k}(\eta_n)
=
\sup_{|I|=k}
\Pp\{|R_{n,I}|>\eta_{n}\}
\le C_k e^{-c_k n^{1/3}}.
\end{equation}
In fact, they state that the constants can be chosen universally over a certain range of $k$,
which is needed for their growing-order Bonferroni argument; we discuss it in the end of this
section.

Next we compute $H_{n,k}(y)$. 
We show that for each fixed $k\ge 1$ and $y\in\R$,
\begin{equation}\label{eq:lixue-H-new}
H_{n,k}(y)\longrightarrow\frac{\{\lambda(y)\}^k}{k!},
\qquad
\lambda(y)=\frac{e^{-y/2}}{\sqrt{8\pi}}.
\end{equation}

\noindent\underline{\underline{\itshape Compute $q_n(y)$ and $H_{n,1}(y)$.}}\hspace{.75em}
Recall
\(
u_n(y)=4\log p_n-\log\log p_n+y\)
and define
\(
q_n(y)=\mathbb P\{V_{1}>u_n(y)\}
\).
We have
\[
q_n(y) = \Pp\left\{\frac{1}{\sqrt{n}}\sum_{i=1}^nX_{i1} X_{i2} >  \sqrt{u_n(y)}\right\} 
       + \Pp\left\{\frac{1}{\sqrt{n}}\sum_{i=1}^nX_{i1} X_{i2} < -\sqrt{u_n(y)}\right\}.
\]
By Cramér-type moderate deviation,
\[
\Pp\left\{\frac{1}{\sqrt{n}}\sum_{i=1}^nX_{i1} X_{i2} > z\right\}=(1+o(1))\overline\Phi(z),\qquad
\Pp\left\{\frac{1}{\sqrt{n}}\sum_{i=1}^nX_{i1} X_{i2} <-z\right\}=(1+o(1))\overline\Phi(z),
\]
uniformly over $0\le z =o(n^{1/6})$.
Taking $z = \sqrt{u_n(y)}$, Mills' ratio yields
\[
2\overline\Phi\!\left(\sqrt{u_n(y)}\right)
\sim \sqrt{\frac{2}{\pi}}\,
\frac{e^{-u_n(y)/2}}{\sqrt{u_n(y)}}
\sim
\sqrt{\frac{2}{\pi}}\,
\frac{p_n^{-2}(\log p_n)^{1/2}e^{-y/2}}{2(\log p_n)^{1/2}}
=\frac{e^{-y/2}}{\sqrt{2\pi}}\,p_n^{-2},
\]
and thus
\begin{equation}\label{eq:lixue-H-n1}
H_{n,1}(y)
=\sum_{\ell}\Pp\{\cE_{n,\ell}(y)\}
=\sum_{\ell}\Pp\{V_{\ell}>u_n(y)\}
=\binom {p_n}2 q_n(y)\longrightarrow
\lambda(y):=\frac{e^{-y/2}}{\sqrt{8\pi}}.
\end{equation}

\noindent\underline{\underline{\itshape Compute $H_{n,k}^{\mathrm{disjoint}}(y)$.}}\hspace{.15em}
We first compute the total contribution from all $k$-sets $I$ for which $\{j_\ell,h_\ell\!:\!\ell\in I\}$ \linebreak contains exactly $2k$ distinct coordinates.
The number of such $k$-sets $I$ is
\[
\binom {p_n}{2k}\frac{(2k)!}{2^k k!}
=\frac{(p_n)_{2k}}{2^k k!},
\]
recalling \((a)_{k}=a(a-1)\cdots(a-k+1)\).
For every such set $I$, we have
\(
\mathbb P\{\mathcal E_{n,I}(y)\}=q_n(y)^k
\)
and thus
\begin{equation}\label{eq:lixue-H-nk-disjoint}
H_{n,k}^{\mathrm{disjoint}}(y)
=\frac{(p_n)_{2k}}{2^k k!}\,q_n(y)^k
=\frac{p_n^{2k}}{2^k k!}
\left(\frac{e^{-y/2}}{\sqrt{2\pi}\,p_n^2}\right)^k+o(1)
=\frac{1}{k!}
\left(\frac{e^{-y/2}}{\sqrt{8\pi}}\right)^k+o(1).
\end{equation}

\noindent\underline{\underline{\itshape Bound $H_{n,k}^{\mathrm{overlap}}(y)$ and compute $H_{n,k}(y)$.}}\hspace{.75em}
The number of $k$-sets $I$ for which $\{j_\ell,h_\ell: \ell\in I\}$ has $2k-1$ or less distinct coordinates is \(O_k(p_n^{2k-1})\).

Consider every such $I$ with $|I|=k$. For any fixed \(\bm\varepsilon=(\varepsilon_{\ell})_{\ell\in I}\in\{-1,1\}^k\),
define
\[
W_{\bm\varepsilon, i}=\sum_{\ell\in I}\varepsilon_{\ell}X_{i j_\ell}X_{i h_\ell}.
\]
Denoting \(S_{\ell}=\frac{1}{\sqrt n}\sum_{i=1}^n X_{ij_\ell} X_{ih_\ell}\), notice that
\[
\varepsilon_{\ell}S_{\ell}>\sqrt{u_n(y)}\quad\text{for all\ \ }\ell\in I
\qquad\Longrightarrow\qquad
\frac1{\sqrt n}\sum_{i=1}^n W_{\bm\varepsilon, i}>k\sqrt{u_n(y)} .
\]
\newpage
\noindent We also have $\Ee W_{\bm\varepsilon, i}=0$ and $\operatorname{Var}(W_{\bm\varepsilon, i})=k$ since $\{j_\ell, h_\ell\} \ne \{j_{\ell'}, h_{\ell'}\}$ for $\ell\ne \ell'$,
and thus by \eqref{eq:pearson-assumptions}
\[
\log\mathbb E e^{tW_{\bm\varepsilon, i}}
=\frac{k t^2}{2}+O_k(|t|^3).
\]
Summing over \(2^k\) choices of signs and using the Chernoff bound gives
\[
\begin{aligned}
\Pp\{\mathcal E_{n,I}(y)\}
=\Pp\Bigg\{\bigcap_{\ell\in I}\Big\{V_{\ell}>u_n(y)\Big\}\Bigg\}
&=\sum_{\bm\varepsilon \in\{-1,1\}^k}\Pp\Bigg\{\bigcap_{\ell\in I}\Big\{\varepsilon_{\ell}S_{\ell}>\sqrt{u_n(y)}\Big\}\Bigg\}\\
&\le
2^k\exp\left\{
-\frac{k u_n(y)}2
+O_k\!\left(\frac{u_n(y)^{3/2}}{\sqrt n}\right)
\right\}
 \le C_{k,y}\,p_n^{-2k}(\log p_n)^{k/2}.
\end{aligned}
\]
This together with the fact that the number of such $I$ is \(O_k(p_n^{2k-1})\) yields
\begin{equation}\label{eq:lixue-H-nk-overlap}
H_{n,k}^{\mathrm{overlap}}(y)
\le C_{k,y}\frac{(\log p_n)^{k/2}}{p_n}
=o(1).
\end{equation}
Combining \eqref{eq:lixue-H-nk-disjoint} and \eqref{eq:lixue-H-nk-overlap} concludes
\[
H_{n,k}(y)
=
\frac{(p_n)_{2k}}{2^k k!}\,q_n(y)^k
+
O_{k,y}\!\left(\frac{(\log p_n)^{k/2}}{p_n}\right).
\]

The remaining ingredient of Corollary~\ref{cor:uniformproxy} is the Lévy concentration probability
mass. Lemma~1 of \citet*{LiXue2015} gives a normal approximation of $F_{n}$ whose Kolmogorov error
is of order $\eta_{n}$. Therefore
\begin{equation}\label{eq:lixue-omega}
\omega_{n}(2\eta_{n})
=
\sup_{x\in\R}\{F_{n}(x+2\eta_{n})-F_{n}(x-2\eta_{n})\}
=O(\eta_{n}),
\end{equation}
because the standard normal distribution has a bounded density;
this is stronger than necessary---the convergence \(F_n\to\Phi\) alone suffices to show that
\(\omega_n(2\eta_n)=o(1)\).

Putting \eqref{eq:lixue-rho}, \eqref{eq:lixue-H-new}, \eqref{eq:lixue-omega} into
Corollary~\ref{cor:uniformproxy} gives the bound
\begin{align}
D_{n,k}(y)
&\le
\rho_{n,k}(\eta_n)\,H_{n,k}(y)
+\rho_{n,k}(\eta_n)\,\binom {d_n}k
+\omega_n(2\eta_n)\times H_{n,k}(y)\notag\\
&=o(1)
\qquad\text{for every fixed }k.
\notag
\end{align}
The asymptotic independence follows by noticing that \eqref{eq:lixue-H-new} also implies the high-order
condition \eqref{eq:H-tail} in Theorem~\ref{thm:bonferroni}.
The joint limit law follows from Corollary~\ref{cor:poisson} by plugging \eqref{eq:lixue-H-new} into
\eqref{eq:poisson-fm} and combining with the marginal limit law for $T_n$.

Most existing literature on high-dimensional independence testing based on max-type statistics
(e.g., \citet*{Jiang2004,Zhou2007,HanChenLiu2017,DrtonEtAl2020})
establishes the limiting distributions of the maxima by deriving the limits of $H_{n,1}$, 
together with Poisson approximation via the Chen--Stein method \citep*{ArratiaEtAl1989,ArratiaEtAl1990}.
Here, by directly calculating $H_{n,k}(y)$ for every fixed $k$, we also provide an alternative way
to establish the limiting distribution of the maximum through Corollary~\ref{cor:poisson}.

The fixed-order Bonferroni argument above is intended only to establish the joint limit,
whereas the growing-order Bonferroni argument in \citet*{LiXue2015} is still needed to establish the
uniform convergence rate. Our framework also accommodates the latter.

\newpage
\subsection{High-dimensional independence testing: finite-$L_q$ blocks of pairwise rank-based statistics}
\label{subsec:zhaowangfeng}

Theorem~2.2 of \citet*{ZhaoWangFeng2026} gives a natural multivariate application of the
local-deletion principle. They also aim to test mutual independence in high dimensions,
but rely on pairwise rank statistics. We retain the notation from the previous section:
write the $p_n$ data columns as $\bm X_1,\ldots,\bm X_{p_n}$ and index the $d_n=\binom {p_n}2$ pairs
in any ordering by
\[
(j_\ell,h_\ell),\qquad 1\le\ell\le d_n.
\]
Let $U_{\ell}$ denote a generic rank-based statistic, which is either
\begin{enumerate}[itemsep=-.1em, topsep=.5em]
\item[(V)] a simple linear rank statistic (e.g., Spearman's~$\rho$ \citep*{Spearman1904}), or 
\item[(W)] a non-degenerate rank-based U-statistic (e.g., Kendall's~$\tau$ \citep*{Kendall1938}), or
\item[(Q)] a degenerate rank-based U-statistic (e.g., Hoeffding's~$D$ \citep*{Hoeffding1948}, \\
    Blum--Kiefer--Rosenblatt's~$R$ \citep*{BlumKieferRosenblatt1961}, \\
    Bergsma--Dassios--Yanagimoto's $\tau^*$ \citep*{BergsmaDassios2014,Yanagimoto1970}),
\end{enumerate}
computed from $\bm X_{j_\ell}$ and $\bm X_{h_\ell}$ in each case.
Let $\widetilde U_{\ell}=U_{\ell}-\Ee(U_{\ell})$;
here and throughout, expectations and probabilities are taken under the null.
Since $U_{\ell}$ is rank-based, we can assume without loss of generality that
all the marginals in the underlying data-generating process $\{X_{ij}:1\le i\le n,\ 1\le j\le p_n\}$
are \linebreak $\text{Unif}\,[0,1]$, and, in particular, are i.i.d. $\text{Unif}\,[0,1]$ under the null.

They define, for each fixed even $q$,
\vspace{-2.5pt}\[\vspace{-2.5pt}
S_{n;q}^{U}
=
\sum_{\ell=1}^{d_n}
\{\widetilde U_{\ell}^q-\Ee(\widetilde U_{\ell}^q)\},
\qquad
\sigma_{n;q}^{U}
=
\sqrt{\sum_{\ell=1}^{d_n}\!\operatorname{Var}(\widetilde U_{\ell}^q)},
\qquad
T_{n;q}^{U}
=
\frac{S_{n;q}^{U}}{\sigma_{n;q}^{U}}.
\]
For a fixed finite set $\cQ=\{q_1,\ldots,q_K\}\subseteq2\mathbb Z_+$,
the sum-type object is defined as the vector
\[
\bm T_{n;\cQ}^{U}
=
(T_{n;q_1}^{U},\ldots,T_{n;q_K}^{U})^\top.
\]

For types (V) and (W), define their max statistics as
\vspace{-2.5pt}\[\vspace{-2.5pt}
M_n^U=\frac{\max_{1\le\ell\le d_n}\widetilde U_\ell^2}
{\operatorname{Var}(\widetilde U_1)}
\qquad
\text{and set}
\quad
\cE_{n,\ell}^U(y)
=
\left\{
\frac{\widetilde U_\ell^2}{\operatorname{Var}(\widetilde U_1)}
>u_n(y)
\right\}
\]
and \(\cE_{n,I}^{U}(y)=\bigcap_{\ell\in I}\cE_{n,\ell}^{U}(y)\) for the extreme-value threshold $u_n(y)=4\log p_n-\log\log p_n+y$, as in \eqref{eq:squareduny}.
The growing regime is given by \eqref{eq:growing-regime}.

\enlargethispage*{0.75\baselineskip}

When $U$ is of type (Q)---a degenerate rank-based U-statistic with a kernel $h$ of order $m$ that is rank-based, symmetric, and satisfies the following three properties: 
\begin{enumerate}[label=(\roman*),itemsep=-.5ex]
\item $h$ is bounded.
\item $h$ is mean-zero and degenerate under independent continuous margins, i.e.,
  $h_1\big( (x_{11}, x_{12}) \big) = 0$.
\item $h_2\big( (x_{11}, x_{12}), (x_{21}, x_{22}) \big)$ has uniformly bounded eigenfunctions, that is, it admits the expansion
\vspace{-2.5pt}\begin{equation*}\vspace{-2.5pt}
h_2\big( (x_{11}, x_{12}), (x_{21}, x_{22}) \big)
=
\sum_{v=1}^{\infty}\lambda_v
\phi_v\big( (x_{11}, x_{12}) \big)
\phi_v\big( (x_{21}, x_{22}) \big),
\end{equation*}
where $\{\lambda_v\}$ and $\{\phi_v\}$ are the eigenvalues and eigenfunctions satisfying the integral equation
\begin{equation*}
\Ee h_2\big( (x_{11}, x_{12}), (X_{21}, X_{22}) \big)
                     \phi\big( (X_{21}, X_{22}) \big)
=
\lambda\phi\big( (x_{11}, x_{12}) \big)~~~\text{for all } (x_{11}, x_{12})\in\R^2,
\end{equation*}
with $\lambda_1\ge\lambda_2\ge\cdots\ge0$, $\Lambda:=\sum_{v=1}^{\infty}\lambda_v\in(0,\infty)$, and $\sup_v\lVert \phi_v\rVert_\infty<\infty$.
\end{enumerate}
Here we write $h_1\big( (x_{11}, x_{12}) \big)                   
  := \Ee h\big( (x_{11}, x_{12}), (X_{21}, X_{22}), \dots, (X_{m1}, X_{m2}) \big)$ \\
  and $h_2\big( (x_{11}, x_{12}), (x_{21}, x_{22}) \big) 
  := \Ee h\big( (x_{11}, x_{12}), (x_{21}, x_{22}), (X_{31}, X_{32}) \dots, (X_{m1}, X_{m2}) \big)$; recall expectations and probabilities are always taken under the null.
The max statistic for type (Q) is accordingly defined as \looseness=-1
\[
M_n^U
=
\frac{n-1}{\lambda_1\binom{m}{2}}
\max_{1\le\ell\le d_n}\widetilde U_\ell
\qquad
\text{and set}
\quad
\cE_{n,\ell}^U(y)
=
\left\{
\frac{n-1}{\lambda_1\binom{m}{2}}\widetilde U_\ell
>u_n^{U}(y)
\right\}
\]
and \(\cE_{n,I}^{U}(y)=\bigcap_{\ell\in I}\cE_{n,\ell}^{U}(y)\) for the extreme-value threshold
\[
u_n^{U}(y) = 4\log p_n+(\mu_1-2)\log\log p_n-\frac{\Lambda}{\lambda_1}+y,
\]
and the growth regime is given by
\begin{equation}\label{eq:growing-regime-theta}
n,p_n\to\infty,\qquad p_n\ge n,\qquad
\log p_n=o(n^{\theta}),
\end{equation}
where $\theta$ is characterized by Equation~(4.1) in \citet*{DrtonEtAl2020}
and can be taken as $1/8-\delta$ \linebreak for any arbitrarily small, prespecified constant $\delta > 0$
for the aforementioned three statistics $D$, $R$, and~$\tau^*$. \looseness=-1

We also impose the remaining regularity assumptions for types (V), (W), and (Q); we refer the reader to the original paper for details.

We use the multivariate version of Corollary~\ref{cor:uniformproxy} to verify the proxy criterion.
Fix $k\ge1$ and $y\in\R$. For $I\subseteq\llbracket d_n\rrbracket$ with $|I|=k$,
$W_I$ is defined the same as in \citet*{LiXue2015} with $|W_I|\le 2k p_n$.
For $q\in\cQ$, define the proxy and the remainder by
\[
\big(S_{n,I;q}^{U}\big)^{\circ}
=
\sum_{\ell\notin W_I}
\{\widetilde U_{\ell}^q-\Ee(\widetilde U_{\ell}^q)\},
\qquad
R_{n,I;q}^{U}
=
\sum_{\ell\in W_I}
\{\widetilde U_{\ell}^q-\Ee(\widetilde U_{\ell}^q)\},
\]
and form the proxy vector and the remainder vector
\[
\big(\bm T_{n,I;\cQ}^{U}\big)^\circ
=
\left(
\frac{\big(S_{n,I;q_1}^{U}\big)^{\circ}}{\sigma_{n;q_1}^{U}},\ldots,
\frac{\big(S_{n,I;q_K}^{U}\big)^{\circ}}{\sigma_{n;q_K}^{U}}
\right)^\top,\qquad
\bm R_{n,I;\cQ}^{U}
=
\left(
\frac{R_{n,I;q_1}^{U}}{\sigma_{n;q_1}^{U}},\ldots,
\frac{R_{n,I;q_K}^{U}}{\sigma_{n;q_K}^{U}}
\right)^\top.
\]
Under the null of mutual independence, it holds that
\(
\big(\bm T_{n,I;\cQ}^{U}\big)^{\circ}
\ind
\cE_{n,I}^{U}(y).
\)

We start with bounding $\rho_{n,k}^{U}$. Combining their assumption that all kernel/score functions
considered are uniformly bounded, along with a rank-moment inequality,
their equation~(C.3) gives for every fixed even $\tau\ge4$ and $\eta\in(0,1)$
\begin{align*}
\rho_{n,k}^{U}(\eta)
&=
\sup_{|I|=k}
\Pp\left\{
\|\bm T_{n;\cQ}^{U}-\big(\bm T_{n,I;\cQ}^{U}\big)^{\circ}\|_\infty>\eta
\right\}
\le
C_{\cQ,\tau}^{U}\,
\eta^{-\tau}k^{2\tau}p_n^{-\tau/2}.\yestag\label{eq:zwf-rho}
\end{align*}
The polynomial bound in \eqref{eq:zwf-rho} is not as tight as the exponential one in
\eqref{eq:lixue-rho}, though still strong enough to absorb the factor $\binom{d_n}{k}$.
Indeed, choosing $\tau=6k$ and $\eta_n=p_n^{-1/12}$ in \eqref{eq:zwf-rho} gives, for fixed $k$,
\[
\rho_{n,k}^{U}(\eta_n)\binom{d_n}{k}
\le
C_{\cQ,k}^{U}\,p_n^{6k/12-3k+2k}
=
O(p_n^{-k/2})
\longrightarrow0.
\]

\newpage
Next we compute $H_{n,k}^{U}(y)$. 
We show that for each fixed $k\ge 1$ and $y\in\R$,
\begin{equation}\label{eq:zwf-H-new}
H_{n,k}^{U}(y)\longrightarrow \frac{\{\lambda^{U}(y)\}^k}{k!},
\qquad
\lambda^{U}(y):=
\begin{cases}
\displaystyle \frac{e^{-y/2}}{\sqrt{8\pi}},
&U\text{ is of type (V) or (W)},\\[8pt]
\displaystyle
\frac{2^{\mu_1/2-2}\kappa}{\Gamma(\mu_1/2)}e^{-y/2},
&U\text{ is of type (Q)},
\end{cases}
\end{equation}
where \(\mu_1:=\#\{v:\lambda_v=\lambda_1\}\)
and \(\kappa:=\prod_{v>\mu_1}\left(1-{\lambda_v}/{\lambda_1}\right)^{-1/2}\).

\noindent\underline{\underline{\itshape Compute $q_n^{U}(y)$ and $H_{n,1}^{U}(y)$.}}\hspace{.75em}
Define
\(
q_n^{U}(y)=\mathbb P\{\cE_{n,1}^{U}(y)\}
\).
By Cramér-type moderate deviation arguments in \citet*{HanChenLiu2017} and 
Cramér-type moderate deviation bounds for degenerate U-statistics in \citet*{DrtonEtAl2020}, respectively,
we obtain
\begin{equation}
q_n^{U}(y)\sim
\begin{cases}
\displaystyle \frac{e^{-y/2}}{\sqrt{2\pi}}\,p_n^{-2},
&U\text{ is of type (V) or (W)},\\[8pt]
\displaystyle
\frac{2^{\mu_1/2-1}\kappa}{\Gamma(\mu_1/2)}e^{-y/2}\,p_n^{-2},
&U\text{ is of type (Q)}.
\end{cases}
\end{equation}
Thus
\begin{equation*}
H_{n,1}^{U}(y)
=\sum_{\ell}\Pp\{\cE_{n,\ell}^{U}(y)\}
=\binom {p_n}2 q_n^{U}(y)\longrightarrow
\lambda^{U}(y).
\end{equation*}

\noindent\underline{\underline{\itshape Compute $(H_{n,k}^{U})^\mathrm{disjoint}(y)$.}}\hspace{.75em}
We first compute the total contribution from all $k$-sets $I$ such that $\{j_\ell,h_\ell: \ell\in I\}$ have $2k$ distinct elements.
The number of such $k$-sets $I$ is
\[
\binom {p_n}{2k}\frac{(2k)!}{2^k k!}
=\frac{(p_n)_{2k}}{2^k k!},
\]
recalling \((a)_{k}=a(a-1)\cdots(a-k+1)\).
For every such set $I$, we have
\(
\mathbb P\{\cE_{n,I}^{U}(y)\}=q_n^{U}(y)^k
\)
and thus \looseness=-1
\begin{equation}\label{eq:zwf-H-nk-disjoint}
(H_{n,k}^{U})^\mathrm{disjoint}(y)
=\frac{(p_n)_{2k}}{2^k k!}\,q_n^{U}(y)^k
=\frac{1}{k!}\lambda^{U}(y)^k+o(1).
\end{equation}

\noindent\underline{\underline{\itshape Bound $(H_{n,k}^{U})^\mathrm{overlap}(y)$ and compute $H_{n,k}^{U}(y)$.}}\hspace{.55em}
We next consider the $k$-sets $I$ for which $\{j_\ell,h_\ell:\ell\in I\}$ \linebreak contains $v\le 2k-1$ distinct coordinates.
We consider the graph on vertex set $\llbracket p_n\rrbracket$ with edges given by the pairs $(j_\ell,h_\ell)$, $\ell\in I$.

Rank-based statistics have many desirable properties concerning independence (see, e.g., Lemma~2.1 in \citet*{LeungDrton2018}),
and using those stated in Proposition~2.1 of \citet*{DrtonEtAl2020} gives the following ``\emph{forest independence}'':
\begin{equation*}
\mathbb P\left(\bigcap_{(j_\ell,h_\ell)\in E(\mathcal F)}\cE_{n,\ell}^{U}(y)\right)
=q_n^{U}(y)^{|E(\mathcal F)|}.
\end{equation*}
This follows by successively removing leaves from a forest \(\mathcal F\), 
conditioning on the remaining variables at each step.
This ``forest independence'' says, as long as the underlying graph induced by $I$ is acyclic, the joint probability of the exceedance events indexed by $I$ factors into the product of their individual probabilities, even when the pairs $(j_\ell,h_\ell)$, $\ell\in I$, share coordinates.

\newpage
Now consider the graph associated with a $k$-set $I$.
Let $v$ and $c$ denote the numbers of non-isolated vertices and nontrivial connected components, respectively.
A spanning forest $\mathcal F_I$ of this graph has exactly $v-c$ edges.
For each fixed $v$, we have
at most
\begin{equation}\label{eq:zwf-overlap-counts}
\binom {p_n}v\binom{\binom v2}{k}=O_k(p_n^v)
\end{equation}
such graphs.
Using the ``forest independence'' 
yields
\[
\Pp\left(\bigcap_{\ell\in I}\cE_{n,\ell}^{U}(y)\right)
\le \Pp\left(\bigcap_{(j_\ell,h_\ell)\in E(\mathcal F_I)}\cE_{n,\ell}^{U}(y)\right)
= q_n^{U}(y)^{v-c}
=O_{y;k}^{U}\!\left(p_n^{-2(v-c)}\right).
\]
Furthermore, we have $v\ge2c+1$ because
every connected component has at least two vertices, 
and at least one of them has three or more---coordinates in $\{j_\ell,h_\ell:\ell\in I\}$ are not distinct.
Therefore
\begin{equation}\label{eq:zwf-overlap-prob}
\Pp\left(\bigcap_{\ell\in I}\cE_{n,\ell}^{U}(y)\right)
=O_{y;k}^{U}\!\left(p_n^{-v-1}\right).
\end{equation}
Taking \eqref{eq:zwf-overlap-counts} and \eqref{eq:zwf-overlap-prob} together yields 
\begin{equation}\label{eq:zwf-H-nk-overlap}
(H_{n,k}^{U})^\mathrm{overlap}(y)
\le \sum_{v=3}^{2k-1}
O_k(p_n^v)\,O_{y;k}^{U}(p_n^{-v-1})
=O_{y;k}^{U}(p_n^{-1})=o(1).
\end{equation}
Hence, combining \eqref{eq:zwf-H-nk-disjoint} and \eqref{eq:zwf-H-nk-overlap},
we finish the proof of \eqref{eq:zwf-H-new}.

We are left with controlling the anti-concentration term.
Let $F_n^{U}(\bm x):=\Pp\{\bm T_{n;\cQ}^{U}\le\bm x\}$,
with $F_{n,j}^{U}(x)$ denoting its $j$th marginal,
and let $\omega_n^U$ denote the corresponding multivariate concentration function.
The marginal Gaussian approximations in Theorem~2.1 of \citet*{ZhaoWangFeng2026} gives
\begin{equation}\label{eq:zwf-omega}
\omega_n^U(2\eta_n)
\le
2\sum_{j=1}^K
\sup_{x\in\R}|F_{n,j}^U(x)-\Phi(x)|
+\frac{4K\eta_n}{\sqrt{2\pi}}
=o(1).
\end{equation}
Combining \eqref{eq:zwf-rho}, \eqref{eq:zwf-H-new}, and \eqref{eq:zwf-omega} gives, for every fixed $k$,
\begin{align*}
D_{n,k}^{U}(y)
&\le
\rho_{n,k}^{U}(\eta_n)\,H_{n,k}^{U}(y)
+\rho_{n,k}^{U}(\eta_n)\,\binom {d_n}k
+\omega_n^{U}(2\eta_n)\times H_{n,k}^{U}(y)
=o(1)
\end{align*}
Since \eqref{eq:zwf-H-new} implies \eqref{eq:H-tail},
Theorem~\ref{thm:bonferroni} yields asymptotic independence.
Whenever the block $\bm T_{n;\cQ}^{U}$ converges in distribution, the joint limit law follows from the multivariate version of Corollary~\ref{cor:poisson}, since \eqref{eq:zwf-H-new} verifies \eqref{eq:poisson-fm}.

\enlargethispage*{.75\baselineskip}

\section{Discussion and extensions}

Theorem~\ref{thm:bonferroni} and Theorem~\ref{thm:proxy} separate a sum--max asymptotic-independence
proof into a generic combinatorial layer and a model-specific approximation layer.
The same philosophy should extend from maxima to extreme point processes.
Instead of a single family of threshold events $\cE_{n,i}(y)$, one may work with factorial moment
measures for counts in finitely many disjoint extreme-value sets.
An independent-proxy criterion of the same type should then yield asymptotic independence between
the sum-type statistic and the entire limiting point process.
Such a result would subsume sum--max asymptotic independence by the continuous mapping theorem
and would connect naturally with the point-process approach of \citet*{HeinyKleemann2023}.
We leave this extension for future work.

\newpage
\section*{Acknowledgment}
The author acknowledges the use of ChatGPT (OpenAI) for language editing and manuscript refinement.

\end{document}